\documentclass{book}
\usepackage[contrib,lang=british]{ems-book} 

\usepackage{euscript}
\usepackage{tikz, tikz-cd, tikz-3dplot}
\usepackage{stmaryrd}

\numberwithin{equation}{section}

\usepackage{upgreek}

\newtheorem{thm}{Theorem}[section]
\newtheorem{prop}[thm]{Proposition}
\newtheorem{lemma}[thm]{Lemma}
\newtheorem{cor}[thm]{Corollary}

\theoremstyle{definition} 
\newtheorem{defin}[thm]{Definition}

\newtheorem{example}[thm]{Example}
\newtheorem{setup}[thm]{Setup}

\newtheorem{remark}[thm]{Remark}

\newtheorem{summary}[thm]{Summary}

\newcommand{\Spec}{\operatorname{Spec}}
\newcommand{\Supp}{\operatorname{Supp}}
\newcommand{\Hom}{\operatorname{Hom}}

\newcommand{\RHom}{\operatorname{{\mathbf R}Hom}}
\newcommand{\End}{\operatorname{End}}
\newcommand{\Pic}{\operatorname{Pic}}
\newcommand{\Aut}{\operatorname{Aut}}
\newcommand{\Auteq}{\operatorname{Auteq}}

\newcommand{\coh}{\operatorname{coh}}

\newcommand{\fmod}{\operatorname{mod}}

\newcommand{\proj}{\operatorname{proj}}

\newcommand{\Kb}{\operatorname{K^b}}

\newcommand{\CM}{\operatorname{CM}}

\newcommand{\SL}{\operatorname{SL}}
\newcommand\Db{\mathop{\mathrm{D}^{\mathrm{b}}}}

\newcommand{\LatticePoint}[1][]{%
\begin{tikzpicture}
\filldraw (0,0) circle (2pt);
\end{tikzpicture}
}

\newcommand{\cStab}[1]{\mathrm{Stab}_{#1}^{\kern -0.5pt \circ}\kern -0.2pt}
\newcommand{\nStab}[1]{\mathrm{Stab}_{#1}\kern -0.1pt}
\newcommand{\TitsR}{{\sf Tits}_{\kern 1pt \bR}}
\newcommand{\CTitsR}{{\sf CTits}_{\kern 1pt \bR}}
\newcommand{\CTitsC}{{\sf CTits}_{\kern 1pt \bC}}

\newcommand{\cAut}[1]{\mathrm{Aut}_{#1}^{\kern -0.5pt \circ}\kern -0.2pt}

\newcommand{\llsq}{[\kern -2.5pt [}
\newcommand{\rrsq}{]\kern -2.5pt ]}

\newcommand{\cH}{\mathcal{H}}

\newcommand{\cS}{\mathcal{S}}

\newcommand{\bC}{\mathbb{C}}

\newcommand{\bP}{\mathbb{P}}

\newcommand{\bR}{\mathbb{R}}

\newcommand{\bZ}{\mathbb{Z}}

\newcommand{\h}{\mathfrak{h}}

\renewcommand{\H}{\mathrm{H}}

\newcommand{\scrA}{\EuScript{A}}

\newcommand{\scrC}{\EuScript{C}}
\newcommand{\scrD}{\EuScript{D}}

\newcommand{\scrF}{\EuScript{F}}

\newcommand{\scrH}{\EuScript{H}}
\newcommand{\scrI}{\EuScript{I}}
\newcommand{\scrJ}{\EuScript{J}}

\newcommand{\scrN}{\EuScript{N}}
\newcommand{\scrO}{\EuScript{O}}

\newcommand{\scrR}{\EuScript{R}}
\newcommand{\scrS}{\EuScript{S}}
\newcommand{\scrT}{\EuScript{T}}

\newcommand{\scrX}{\EuScript{X}}
\newcommand{\scrY}{\EuScript{Y}}
\newcommand{\scrZ}{\EuScript{Z}}
\newcommand{\con}{\operatorname{con}}

\newcommand{\scrIc}{\scrI^{\kern 0.5pt\mathrm{c}}}

\newcommand{\Delt}{\Updelta_{\kern 0.05em 0}}
\newcommand{\DeltAff}{\Updelta_{\kern 0.05em 0}^{\aff}}
\newcommand{\WDelt}{W_{\kern -0.1em \Updelta}}
\newcommand{\WDeltaff}{W_{\kern -0.1em \Updelta_{\aff}}}
\newcommand{\Wkern}[1]{W_{\kern -0.1em #1}\kern 0.05em}
\newcommand{\wo}[1]{w_{\kern -0.075em #1}}
\newcommand{\wop}[1]{w^{\phantom J}_{\kern -0.1em #1}}
\newcommand{\GammaJ}{\Upgamma_{\kern -0.05em J}}
\newcommand{\GammaS}{\Upgamma_{\kern -0.05em \scrJ}}

\newcommand{\xlup}[1]{{}^{#1}\kern -0.15em x}
\newcommand{\xlupmax}[1]{{}^{#1}\kern -0.25em x}
\newcommand{\iDelta}{\iota_{\kern -0.075em \Updelta}}
\newcommand{\Phisub}[1]{\Phi_{\kern -0.1em #1}}
\newcommand{\aff}{\operatorname{\mathtt{af}}\nolimits}

\definecolor{Pink}{RGB}{230 56 243}
\tikzset{
        pvertex/.style={circle,inner sep=1pt,outer sep=2pt,font=\scriptsize},
W/.style={circle,draw=black,circle,fill=white,inner sep=0pt, minimum size=4pt},
B/.style={circle,draw=black!80!white,circle,fill=black!80!white,inner sep=0pt, minimum size=4pt},
Or/.style={circle,draw=Orange,circle,fill=Orange,inner sep=0pt, minimum size=4pt},
P/.style={circle,draw=Pink,circle,fill=Pink,inner sep=0pt, minimum size=4pt},
R/.style={circle,draw=black!80!white,circle,fill=red!80!white,inner sep=0pt, minimum size=4pt},  
}

\newcommand{\Esix}[6]{%
\begin{tikzpicture}[scale=0.21]
\node at (0,0) [#1] {};
\node at (1,0) [#2] {};
\node at (2,0) [#3] {};
\node at (2,1) [#4] {};
\node at (3,0) [#5] {};
\node at (4,0) [#6] {};
\end{tikzpicture}
}
\newcommand\Curve{\mathtt{C}}

\def\Id{\mathop{\mathrm{Id}}\nolimits}

\newcommand\chamC{\mathrm{C}}
\newcommand\chamD{\mathrm{D}}

\newcommand{\smc}{\operatorname{\mathrm{smc}}}

\newcommand{\silt}{\operatorname{\mathrm{silt}}}

\newcommand{\Hyp}{\mathsf{H}}

\begin{document}
\mainmatter

\title{Spherical and Semibrick Classifications}
\titlemark{Spherical and Semibrick Classifications}

\emsauthor{1}{
	\givenname{Wahei}
	\surname{Hara}
	\mrid{1228979}
	\orcid{0009-0003-6567-8496}}{W.~Hara}
\emsauthor{2}{
	\givenname{Michael}
	\surname{Wemyss}
	\mrid{893224}
     \orcid{0000-0003-4760-7750}}{M.~Wemyss} 

\Emsaffil{1}{
	\pretext{}
	\department{Kavli Institute for the Physics and Mathematics of the Universe (WPI)}
	\organisation{University of Tokyo}
	\rorid{01a2bcd34}
	\address{5-1-5 Kashiwanoha, Kashiwa}
	\zip{277-8583}
	\city{Chiba}
	\country{Japan}
	\posttext{}
	\affemail{wahei.hara@ipmu.jp}
	\furtheremail{}}
 \Emsaffil{2}{
	\department{School of Mathematics and Statistics}
	\organisation{University of Glasgow}
	\rorid{00vtgdb53} 
	\address{University Place}
	\zip{G12 8QQ}
	\city{Glasgow}
	\country{UK}
	\affemail{michael.wemyss@glasgow.ac.uk}}

\classification[14A22, 14E16, 14J30, 16D60]{18G80}



\begin{abstract}

This article provides an overview of the techniques related to classification of spherical and more general objects within triangulated categories, and its relationship with algebraic geometry, representation theory and symplectic geometry. The primary focus are the techniques within the geometric but `finite' setting of \cite{HaraWemyss}, but other approaches including to more general settings by \cite{IU, BDL, SW, KeatingSmith, Parth2} are also surveyed, in varying levels of detail. Various explicit examples are provided. 
 
\end{abstract}

\makecontribtitle



\begin{ack}
We thank Tom Bridgeland, Osamu Iyama, Ailsa Keating, Tony Licata, David Pauksztello, Parth Shimpi and Ivan Smith for conversations and collaborations, the conference organisers for inviting us to submit to this ICRA volume, and the referee for careful reading.
\end{ack}

\begin{funding}
The authors were supported by EP/R034826/1, and by ERC Consolidator Grant 101001227 (MMiMMa).
W.H.~was further supported by World Premier International Research Center Initiative (WPI), MEXT, Japan, and by JSPS KAKENHI Grant Number JP24K22829. \end{funding}

\section{Introduction}

Through time, various phenomena in mathematics become naturally linked. Motivated by ideas and approaches to Brou\'e{}'s conjecture in modular representation theory in the 1990s, particularly Okuyama’s method, in \cite[5.1]{Rickard} Rickard characterised `simple-like' objects in the derived category of a finite dimensional symmetric algebra.  The key observation was that if some collection of generating objects looks `a bit like simples', they are simples, after passing through some shiny new derived equivalence.  A follow-up paper of \cite{AlN} refined these ideas further, generalising in two directions: first by dropping the symmetric assumption, then second by replacing the derived equivalence with a more general statement on t-structures. 

These generating collections of `simple-like' objects were called \emph{cohomologically Schurian sets of generators}, and perhaps mercifully, through time this language has softened to the linguistically easier \emph{simple-minded collections}. From the viewpoint of what follows, we pause here to dwell slightly on the loss of the word `generator' from the lexicon, since this key technical point is perhaps the singlemost distinguishing factor that is lacking in all the geometric phenomena below. 

\medskip
Around the same time, but instead motivated by mirror symmetry considerations, spherical objects entered the fray \cite{ST}.  For a given $k$-linear triangulated category $\scrT$, these are by definition objects $x\in\scrT$ which share the same cohomology as a sphere, in the sense that
\[
\Hom_{\scrT}(x,x[t])\cong\begin{cases}k&\mbox{if }t=0,d\\ 0&\mbox{else,}\end{cases}
\]
together with some additional compatibility with the Serre functor.  Modulo the usual caveats regarding suitable enhancements, the key property of spherical objects $x$ is that they generate an equivalence $T_S\colon\scrT\to\scrT$, called the spherical twist.

Delaying slightly the big reveal, there are multiple reasons why spherical objects are important, many of which we will not list here. The classification of spherical objects became a desirable theme, partly since in the setting of Fukaya categories various remarkable topological consequences immediately follow. Similarly, in various algebraic-geometric settings there are immediate applications to the structure of the autoequivalence group, and to the connectedness of the stability manifold.

\medskip
By definition spherical objects are `a bit like a sphere', however as observed by many authors, they too have the habit of behaving `a bit like a simple'.  From this, the obvious analogy is born.  Much like the very best analogies, from the outset this is completely ridiculous: spherical objects are \emph{single} objects, and whilst Rickard demonstrated that having a generating set of simples is desirable, it was well-known that having a single simple object is generally as good as useless.

\medskip
The above point makes the need to restrict the setting obvious.  Below we survey various categories in which both the classification of single simple objects \emph{and} of generating sets of simple objects is possible, but in the interest of full disclosure, we need to be clear that there is currently no single abstract framework in which all these settings naturally belong. The techniques surveyed do transfer between the settings, after some reinterpretation, but work is required to do this.  


Although the axiomatics currently escapes us, we do contest that all the categories~$\scrC$ (although not the categories $\scrD$) considered in this survey should be viewed as \emph{Iyama dream categories} in a manner that will one day be axiomatically defined.  The main point is that all the categories~$\scrC$ considered below turn out to admit the dream level of completions: in particular, all simple objects can be completed to a simple-minded collection, and so the difference between single simple/spherical objects and generating collections melts away.  The broader reason why the categories $\scrC$ should satisfy this property remains a mystery; currently our techniques only tell us that it is true, they do not give much of a hint of why or when we should expect this behaviour more generally.


\medskip 
The remainder of this article overviews and compares approaches to the classification of spherical objects, in various contexts. Of course, given the above, we in fact survey more general classification results too, and the differing algorithms used in the various cases.  

As one piece of language, many of the categories below will have a  naturally associated hyperplane arrangement.  In the setting that this arrangement is finite, we will denote the category $\scrC$ and refer to it as a `finite' setting, and in the case when the arrangement is infinite, we will denote the category $\scrD$, and refer to it as the `affine' setting.  Given the fact that Iyama dream categories are not yet defined, there are of course some categories considered below (such as $\Db(\fmod A)$ in \S\ref{sec:SiltingDiscrete}) which strictly speaking fall into neither setting, but their simple-like objects can still be understood.

\section{Definitions}
In a triangulated category $\scrT$ the notion of a \emph{brick}, recalled below, is a replacement for the notion of a simple module. Indeed, the vanishing of the negative Ext groups mimicks the behaviour of an object in the heart of a t-structure, whilst the condition that the Hom-space is one-dimensional mimicks Schur's lemma.

\begin{defin} 
Let $\scrT$ be a triangulated category over $\bC$.
\begin{enumerate}
\item $x \in \scrT$ is a \emph{brick}
if $\Hom_\scrT(x,x[i]) = 0$ for all $i < 0$ and $\Hom_\scrT(x,x) = \bC$.
\item A \emph{semibrick} is an object of the form $y = \bigoplus_{i=1}^n x_i$ where each $x_i$ is a brick, and further $\Hom_\scrT(x_i, x_j[k]) = 0$ for all $i \neq j$ and $k \leq 0$.
\item A \emph{simple-minded collection (smc)} in $\scrT$ is a collection of objects $\{x_1, \hdots, x_n \}$ such that $y = \bigoplus_{i=1}^n x_i$ is a semibrick, and further $y$ is a classical generator of $\scrT$, in the sense that $\scrT$ is the smallest triangulated subcategory of $\scrT$ containing $y$ that is closed under taking direct summands.
\end{enumerate}
By abuse of notation, given an smc $\{x_1, \hdots, x_n \}$,
we also call the direct sum $y = \bigoplus_{i=1}^n x_i$ an smc of $\scrT$.
\end{defin}

\begin{defin}\label{def:silting}
Let $\scrT$ be a triangulated category over $\bC$,
and $P \in \scrT$ an object.
\begin{enumerate}
\item[(1)] $P$ is called \emph{presilting} if $\Hom_\scrT(P,P[i]) = 0$ for all $i>0$.
\item[(2)] $P$ is called \emph{silting} if $P$ is presilting, and a classical generator of $\scrT$. 
\end{enumerate}
\end{defin}

\begin{remark}\label{rem:Dual}
For a finite dimensional $\bC$-algebra $A$,
the full collection of simple $A$-modules is an example of a simple-minded collection in $\Db(\fmod A)$, and $A$ itself is an example of a silting object in $\Kb(\proj A)$.
Note that, the notions of smc and silting are in some sense dual to each other, but they exist in different categories.
\end{remark}

\section{Summary of Common Strategy}\label{sec:CommonStrategy}

The aim of this section is to give an overview of the common strategy to most approaches to the classification of spherical (and others) objects. Details appear later, in future sections.  The most common strategy is the following:
\begin{enumerate}
\item Fix a category $\scrT$, and some subset of objects defined by some chosen homological property.
\item Consider some invariant $\ell(x) \in \bR$ of objects $x$ in $\scrT$.
\item Classify all objects $x$ in $\scrT$, which 
\begin{enumerate}
\item have small $\ell(x)$, and 
\item satisfy the homological property in (1).
\end{enumerate}
\item Find an induction step that makes the invariant $\ell(x)$ smaller, by e.g.\ applying functors such as spherical twists, or by varying t-structures.
\end{enumerate}
There are really four choices: the category, the homological property, the invariant, and the induction step. In some sense it has to be taken as given that the classification in (3) is possible, since otherwise the wider problem is hopeless.  Combined, the choices above should be viewed as a choice of an algorithm that has as its outcome the desired classification of objects satisfying the homological condition in (1).


We remark that step (4) can often (but not always, see \ref{rem:notenoughAut}) be thought of as a choice of a set of autoequivalences $\nabla$. From that viewpoint, the existence of an improvement of $\ell(x)$ is equivalent to establishing the non-emptiness of the set
\[ 
\left\{\,\Phi \in \nabla \mid \ell(\Phi x) < \ell(x)   \,\right\}. 
\]
Depending on the choices made in (1)--(3), there is sometimes advantage in having $\nabla$ to be a large set, and sometimes there is advantage in having $\nabla$ to be small. In the works by Ishii-Uehara \cite{IU} and Bapat-Deopurkar-Licata \cite{BDL} surveyed below, $\nabla$ consists of (single) spherical twists along spherical sheaves. On the other hand, in \S\ref{subsec:HomMMP} $\nabla$ is taken to be compositions of such functors, albeit in a controlled way. 

In what follows, we survey the main approaches in settings from algebraic geometry, representation theory and symplectic geometry.  In each case, the four choices above are all subtly different, mainly since they strongly depend on the precise category being considered, and on what precise objects are being classified.

\section{Geometry Settings and HomMMP}\label{HW:mainsec}\label{subsec:HomMMP}
This section is concerned with the following three geometric settings.
\begin{setup}[Geometric setup] \label{geom setup}
Let $f \colon X \to \Spec R$ be one of the following.
\begin{enumerate}
\item\label{geom setup 1} The minimal resolution $f \colon \scrZ \to \bC^2/G$, where $G \leq \SL(2,\bC)$ is finite.
\item\label{geom setup 2} A partial crepant resolution $f \colon \scrY \to \bC^2/G$, where $G \leq \SL(2,\bC)$ is finite.
In other words, $f$ is a crepant birational morphism. 
\item\label{geom setup 3} A threefold flopping contraction $f \colon \scrX \to \Spec \scrR$
such that $\scrR$ is an isolated complete local cDV singularity.
\end{enumerate}
\end{setup}
In (1) $\scrZ$ is smooth, in (2) $\scrY$ has at worst Kleinian singularities, and in~(3) $\scrX$ has at worst isolated cDV singularities. Of course, (1) is really just a special case of (2).

We will consider the \emph{null category}
\[
\scrC\coloneq\{ a\in\Db(\coh X)\mid  \mathbf{R}f_*a=0\},
\]
where $X$ is either $\scrX, \scrY$ or $\scrZ$ above. Later we will also consider the category
\[
\scrD\coloneq\{ a\in\Db(\coh X)\mid \Supp a\subseteq \Curve\},
\]
where $\Curve$ is the exceptional locus, but our primary focus will be on the category $\scrC$. 

\subsection{Hyperplane Arrangements}
For each $X$ as above, we first associate a marked Dynkin diagram $(\Delta, \scrJ)$ as  follows.
\begin{enumerate}
\item For the minimal resolution $\scrZ \to \bC^2/G$, the dual graph of the exceptional curve $\Curve_{\mathrm{red}}$ is a Dynkin diagram $\Delta$ of ADE type. We associate to $\scrZ$ the pair $(\Delta,\emptyset)$.
\item If $\scrY \to \bC^2/G$ is a partial crepant resolution,
then consider the minimal resolution $\scrZ \to \scrY \to \bC^2/G$,
and let $\Delta$ be the ADE Dynkin diagram associated to $\scrZ \to \bC^2/G$.
We associate to $\scrY$ the pair $(\Delta,\scrJ)$, where $\scrJ \subset \Delta$ is the set of vertices that correspond to exceptional curves $\Curve_i \subset \scrZ$ contracted by $\scrZ \to \scrY$.
\item If $\scrX \to \Spec \scrR$ is a threefold flopping contraction,
then for a general element $g \in \scrR$, the quotient $\scrR/g$ is isomorphic to $\bC\llsq  x,y\rrsq ^G$ for some finite $G \leq \SL(2,\bC)$, and the base change $\scrX \otimes_{\scrR} \scrR/(g) \to \Spec \scrR/(g)$ is a partial crepant resolution of the form in (2). Through this, we associate to $\scrX$ the corresponding $(\Delta, \scrJ)$ in (2).
\end{enumerate}

Given the data $(\Delta, \scrJ)$, it is possible to construct a finite hyperplane arrangement $\scrH_\scrJ$. It is classical that $\Delta$ has an associated root system $\h = \bigoplus_{i \in \Delta} \bR \upalpha_i$, where $\upalpha_i$ are the simple roots. Let $\h_{\scrJ} \coloneq {\h}/{\bigoplus_{i \in \scrJ} \bR \upalpha_i}$ be the quotient space corresponding to the marking $\scrJ \subset \Delta$, with the natural projection 
\[ 
\uppi_{\scrJ} \colon \h \to \h_{\scrJ}. 
\]
Note that the space $\h_{\scrJ}$ has basis $\{\uppi_{\scrJ}(\upalpha_j) \mid j \in \Delta\setminus\scrJ \}$.

\begin{defin}
A \emph{restricted positive root} in $\h_{\scrJ}$ is a non-zero element $\uppi_{\scrJ}(\upalpha) \in \h_{\scrJ}$ for some positive root $\upalpha \in \h$.
A restricted positive root is \emph{primitive} if it is not a multiple of another positive restricted root.
\end{defin}

\begin{example}\label{ex:EsixA}
Consider the Dynkin data $\,\Esix{W}{B}{W}{B}{B}{B}$, where $\scrJ$ is the set of black nodes.  Projecting all positive roots of $E_6$ via $\uppi_\scrJ$, and discounting the zeros, gives the set
\[
\{  10, 01, 02, 11, 12, 13 \}
\] 
where for example $11$ is shorthand for the coordinate $(1,1)$ in $\mathfrak{h}_\scrJ$ under the prescribed basis above. These are the restricted positive roots.  The primitive restricted roots are $\{  10, 01, 11, 12, 13 \}$.
\end{example}

\begin{defin}
Let $(\Delta, \scrJ)$ be a marked Dynkin diagram,
and $\h_{\scrJ}$ the associated space.
Set $\Theta_{\scrJ} \coloneq (\h_{\scrJ})^{\vee}$, and for each restricted positive root $\upbeta = \uppi_{\scrJ}(\upalpha)$, consider
\[ 
\Hyp_{\upbeta} \coloneq  \{(\upvartheta_i)\in\Theta_\scrJ \mid \sum \upbeta_i \upvartheta_i = 0 \} \subseteq \Theta_{\scrJ}. 
\]
The set $\scrH_\scrJ \coloneq \{\Hyp_{\upbeta}\mid \upbeta \mbox{ is restricted positive root}\}$ is the \emph{intersection arrangement}.
\end{defin}

\begin{example}\label{ex:EsixB}
Continuing \ref{ex:EsixA}, e.g.\ the restricted root $12$ becomes the hyperplane $\upvartheta_1+2\upvartheta_2=0$.  The full arrangement $\scrH_\scrJ$ is illustrated below.
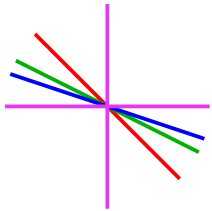
\[
\begin{array}{ccc}
\begin{array}{c}
\begin{tikzpicture}[scale=0.6]
\coordinate (A1) at (135:2.25cm);
\coordinate (A2) at (-45:2.25cm);
\coordinate (B1) at (153.435:2.25cm);
\coordinate (B2) at (-26.565:2.25cm);
\coordinate (C1) at (161.565:2.25cm);
\coordinate (C2) at (-18.435:2.25cm);
\draw[line width=0.5mm, red] (A1) -- (A2);
\draw[line width=0.5mm,green!70!black] (B1) -- (B2);
\draw[line width=0.5mm, blue] (C1) -- (C2);
\draw[line width=0.5mm, Pink] (-2.25,0)--(2.25,0);
\draw[line width=0.5mm, Pink] (0,-2.25)--(0,2.25);
\end{tikzpicture}
\end{array}
&&
\begin{array}{c}
\begin{tabular}{rcc}
\toprule
Hyperplane&&Restricted Roots\\
\midrule
$\upvartheta_1=0$&&$10$\\
$\upvartheta_2=0$&&$01, 02$\\
$\upvartheta_1+\upvartheta_2=0$&&$11$\\
$\upvartheta_1+2\upvartheta_2=0$&&$12$\\
$\upvartheta_1+3\upvartheta_2=0$&&$13$\\
\bottomrule
\end{tabular}
\end{array}
\end{array}
\]
\end{example}

\subsection{Functors}
Let $f \colon X \to \Spec R$ be one of the options in Setup~\ref{geom setup}.
Van den Bergh \cite{VdB} proves the existence of a vector bundle $\scrO\oplus \scrN$ on $X$, and a derived equivalence
\[ 
\Psi_X \colon \Db(\coh X) \xrightarrow{\RHom_X(\scrO \oplus \scrN, -)} \Db(\fmod \Lambda), 
\]
where $\Lambda \coloneq \End_X(\scrO \oplus \scrN)$.  
\begin{remark}
When $\dim X=2$, the noncommutative algebra $\Lambda$ can be described in terms of representation theory.
\begin{enumerate}
\item[(1)] If $X = \scrZ$ is the minimal resolution, then $\Lambda$ is isomorphic to the preprojective algebra $\Uppi$ of the corresponding extended Dynkin diagram.
\item[(2)] If $X = \scrY$ is a partial crepant resolution of a Kleinian singularity,
let $(\Delta, \scrJ)$ be the associated marked Dynkin diagram.
Let $\Uppi$ be as in (1), and $e_0, e_1, \hdots, e_n \in \Uppi$ the idempotents,
where $0$ corresponds to the extended vertex.
Then $\Lambda$ is isomorphic to the contracted affine preprojective algebra $\upvarepsilon_{\scrJ} \Uppi \upvarepsilon_{\scrJ}$, where $\upvarepsilon_{\scrJ} = 1 - \sum_{i \in \scrJ} e_i$.
\end{enumerate}
\end{remark}

In each of the cases of Setup~\ref{geom setup}, we next associate categorical information to the hyperplane arrangement $\scrH_\scrJ$. 

\begin{summary}\label{sum:CategoriesAssignment}
Given a fixed $f \colon X \to \Spec R$ as in Setup~\ref{geom setup}, consider the associated $\scrH_\scrJ$. To this, we associate the following categorical data:
\begin{itemize}
\item to every chamber $\chamC$, a category $\scrC_{\chamC}$,
\item to every wall-crossing $s_i \colon \chamC \to \chamD$, an equivalence $\Phi_i \colon \scrC_{\chamC} \to \scrC_{\chamD}$, 
\end{itemize}
in the following way, depending on the type of $f \colon X \to \Spec R$.

\begin{enumerate}
\item For the minimal resolution $f \colon \scrZ \to \bC^2/G$, with associated $(\Delta,\emptyset)$, the hyperplane arrangement is the classical reflection arrangement of the Weyl group.
To every chamber $\chamC$ we associate the same category, namely
\[ 
\scrC_{\chamC} \coloneq \{ a \in \Db(\coh \scrZ) \mid Rf_*a = 0 \} ,
\]
and to each simple wall crossing $s_i\colon\chamC\to\chamD$ with $i=1,\hdots,n$, is assigned the spherical twist along $\scrO_{\Curve_i}(-1)$,
where $\Curve_i$ is the irreducible component corresponding to the index $i$, which is isomorphic to $\bP^1$.

\item For a partial crepant resolution $\scrY \to \bC^2/G$, with associated $(\Delta,\scrJ)$, we give an algebraic reinterpretation and generalisation of the above. 
Now chambers in $\scrH_\scrJ$ are indexed by certain pairs $(x,\scrI)$ where $x$ belongs to the Weyl group of $\Updelta$, and $\scrI \subset \Updelta$ \cite[1.12]{IW9}.  
For a chamber $\chamC=(x,\scrI)$, set $\Lambda_{\chamC} \coloneq \upvarepsilon_{\scrI} \Uppi \upvarepsilon_{\scrI}$ and associate the category
\[
\scrC_{\chamC} \coloneq \{ a\in\Db(\fmod \upvarepsilon_{\scrI} \Uppi \upvarepsilon_{\scrI})\mid e_0\,\H^*(a)=0 \}.
\]
where $e_0$ is the idempotent corresponding to the extended vertex.  These categories now typically vary, depending on the chamber.  To each simple wall crossing $\upomega_i\colon\chamC\to\chamD$ is an associated mutation functor $\Phi_i$, defined generally in \cite[\S6]{IW1}, and in the setting here in \cite[\S5.6]{IW9}.


\item For a flopping contraction $\scrX\to\Spec\scrR$, by \cite{HomMMP} there is a bijection between chambers of $\scrH_\scrJ$ and the movable cone relative to $f$.  Thus, to every chamber $\chamC$ is an associated $f_\chamC\colon\scrX_\chamC\to\Spec\scrR$, and we thus associate to $\chamC$ the category 
\[
\scrC_{\chamC} \coloneq \{ a \in \Db(\coh \scrX_{\chamC}) \mid R(f_{\chamC})_*a = 0 \}.
\]
To each simple wall crossing $s_i \colon \chamC \to \chamD$,
there is a corresponding curve $\Curve_i$, to which the birational map $\scrX_{\chamC} \dashrightarrow \scrX_{\chamD}$
is the flop of $\Curve_i$. To the wall crossing we thus associate the equivalence given by the (inverse) of the Bridgeland--Chen flop functor \cite{B02, Chen}.

Equivalently, but more algebraically, by \cite{HomMMP, IW9} there is a bijection $\chamC\mapsto M_\chamC$ between chambers in $\scrH_\scrJ$ and certain rigid objects in the category of Cohen--Macaulay modules $\CM\scrR$.  
Necessarily $M_\chamC$ has $\scrR$ as a summand. 
Thus to a chamber $\chamC$ set $\Lambda_{\chamC} \coloneq \End_{\scrR}(M_{\chamC})$, and assign the subcategory 
\[
\scrC_{\chamC} \coloneq \{ a\in\Db(\fmod \Lambda_{\chamC})\mid e_0\,\H^*(a)=0 \}
\]
where $e_0$ is the idempotent corresponding to the summand $\scrR$.
To each simple wall crossing $s_i\colon\chamC\to\chamD$, assign the mutation functor 
\[ 
\Phi_i\coloneq \RHom_{\Lambda_\chamC}(\Hom_\scrR(M_\chamC,M_\chamD),-). 
\]
\end{enumerate}
\end{summary}



The noncommutative algebra $\Lambda_{\chamC}$ with a fixed idempotent $e_0 \in \Lambda_{\chamC}$ as above introduces an important class of finite dimensional algebras, as follows.

\begin{defin}
Given a chamber $\chamC$ in the hyperplane arrangement $\scrH_\scrJ$ associated to $f \colon X \to \Spec R$ in Setup~\ref{geom setup}, consider the noncommutative algebra $\Lambda_{\chamC}$ described in \ref{sum:CategoriesAssignment} (either $\Pi$, $\upvarepsilon_\scrJ\Pi\upvarepsilon_\scrJ$ or $\End_\scrR(M_\chamC)$), and let $e_0 \in \Lambda_{\chamC}$ be the idempotent defined as before. Set
\[ 
\Lambda_{\chamC, \con} \coloneq \Lambda_{\chamC}/(e_0), 
\]
which is known to be finite dimensional over $\bC$.
\end{defin}

This algebra also has a representation-theoretic description.

\begin{enumerate}
\item[(1)] If $X = \scrZ$ is the minimal resolution, then $\Lambda_{\chamC, \con} = \Uppi/(e_0)$ is the preprojective algebra of Type $\Delta$.
\item[(2)] If $X = \scrY$ is a partial crepant resolution, then $\Lambda_{\chamC, \con} = e_{\scrJ} \Uppi e_{\scrJ}/(e_0)$ is the contracted preprojective algebra associated to $(\Delta, \scrJ)$.
\item[(3)] If $X = \scrX \to \Spec \scrR$ is a threefold flopping contraction,
then the algebra $\Lambda_{\chamC, \con}$ is called the \emph{contraction algebra} of $f_{\chamC} \colon \scrX_{\chamC} \to \Spec \scrR$, defined by Donovan-Wemyss \cite{DW1, DW3}.
\end{enumerate}

The algebra $\Lambda_{\chamC, \con}$ provides a bounded $t$-structure 
in the category $\scrC_{\chamC}$. Indeed, the quotient map $\Lambda_{\chamC} \to \Lambda_{\chamC, \con}$
gives an embedding $\fmod \Lambda_{\chamC, \con} \subset \fmod \Lambda_{\chamC}$.
The equality
\[ 
\fmod \Lambda_{\chamC, \con} = {\fmod \Lambda_{\chamC}} \cap {\scrC_{\chamC}} 
\]
holds, and is the heart of bounded $t$-structure in $\scrC_{\chamC}$.

\begin{defin}
The heart $\fmod \Lambda_{\chamC, \con} \subset \scrC_{\chamC}$
is called the \emph{standard heart} of $\scrC_{\chamC}$.
\end{defin}

\begin{remark} \label{rem simples of std heart}
The standard heart has exactly $n$ simple modules $\scrS_1, \hdots, \scrS_n$,
where $n$ is equal to the number of irreducible curves $\Curve_i$. Across the equivalence $\Psi_{X_{\chamC}}$ between $\Db(\coh \scrX_{\chamC})$ and $ \Db(\fmod \Lambda_{\chamC})$, the simple module $\scrS_i$ corresponds to $\scrO_{\Curve_i}(-1)$. 
 In addition, the mutation functor $\Phi_i$ satisfies $\Phi_i(\scrS_i) = \scrS_i[-1]$. 
\end{remark}

\subsection{The Classification Result}\label{subsec:1.4CLass}

With the above in place, we now state the main classification result.

\begin{thm} \label{HaraWemyss main thm}
Let $f \colon X \to \Spec R$ be as in Setup~\ref{geom setup},
$\scrH_\scrJ$ the associated hyperplane arrangement, and $\chamC$ a chamber. Then the following statements hold.
\begin{enumerate}
\item\label{HaraWemyss main thm 1} An object $x \in \scrC_{\chamC}$ satisfies $\Hom(x,x[i]) = 0$ for all $i < 0$ if and only if there is a composition of mutation functors $\Phi_{\upalpha}\colon\scrC_\chamC\to\scrC_\chamD$, for some chamber $\chamD$, such that $\Phi_{\upalpha}(x)$ is contained in the standard heart $\fmod \Lambda_{\chamD, \con} \subset \scrC_{\chamD}$.
\item\label{HaraWemyss main thm 2} An object $x \in \scrC_{\chamC}$ satisfies $\Hom(x,x[i]) = 0$ for all $i < 0$ and $\Hom(x,x) = \bC$ (i.e.~$x$ is a brick) if and only if  there is a composition of mutation functors $\Phi_{\upalpha}\colon\scrC_\chamC\to\scrC_\chamD$, for some chamber $\chamD$, such that $\Phi_{\upalpha}(x)\cong \scrO_{\Curve_j}(-1)$
for some~$j$.
\end{enumerate}
\end{thm}
The proof deploys the common strategy, using the following invariant $\ell(x)$.

\begin{defin}\label{defin:abNotation}
Consider the category $\scrC_{\chamC}$ associated to the chamber $\chamC$, with standard heart $\fmod \Lambda_{\chamC,\con}$, and write $\H^* \colon \scrC_{\chamC} \to \fmod \Lambda_{\chamC,\con}$ for the cohomology functor associated with this t-structure. Then, for an object $x \in \scrC_\chamC$ and integers $a \leq b \in \bZ$,
\begin{enumerate}
\item $x \in [a,b]$ indicates that $\H^i(x) = 0$ for all $i < a$ and $i>b$.
\item $x \in \llsq a,b\rrsq $ indicates that $x \in [a,b]$ and $\H^i(x) \neq 0$ for $i = a, b$.
\end{enumerate}
If an object $x \in \scrC_{\chamC}$ satisfies $x \in \llsq a,b\rrsq $, put
\[ \ell(x) \coloneq b -a. \]
\end{defin}

An object $x \in \scrC_{\chamC}$ satisfies $\ell(x) = 0$
if and only if $x$ belongs to a shift of the standard heart.
An important step for proving Theorem~\ref{HaraWemyss main thm}
is to guarantee that one can make $\ell(x)$ smaller by applying a composition of mutation functors.
To show this, one may want to test $\ell(\Phi_ix)$ among all single mutations $\Phi_i$, however this naive strategy does not work well.

The correct strategy is to test $\ell(\Phi_{\upalpha}x)$ among all \emph{atoms} $\upalpha$.
Atoms can be characterized using the standard hearts as follows,
and can be seen as a natural class of functors that contains the identity and all single mutation functors $\Phi_i$.

\begin{lemma}
Let $f \colon X \to \Spec R$ be as in Setup~\ref{geom setup},
with corresponding hyperplane arrangement $\scrH_\scrJ$. 
\begin{enumerate}
\item A composition of mutation functors $\Phi_{\upalpha}\colon\scrC_\chamC \to \scrC_\chamD$ is an atom if and only if $\Phi_{\upalpha}(y) \in [0,1]$ for all $y \in \fmod \Lambda_{\chamC, \con}$.

\item For any two chambers $\chamC$ and $\chamD$, there exists a unique atom
$\Phi_{\upalpha} \colon \scrC_\chamC \to \scrC_\chamD$.
In particular, since there are finitely many chambers,
the set of all atoms is finite.
\end{enumerate}
\end{lemma}

The set of atoms admits a natural partial order, which plays an important role.

\begin{defin}
Consider two atoms $\Phi_{\upalpha_i}\colon\scrC_{\chamC_i} \to \scrC_{\chamD_i}$ for $i = 1,2$.
Then $\Phi_{\upalpha_1} > \Phi_{\upalpha_2}$ indicates that
\begin{enumerate}
\item $\chamC_1 = \chamC_2$, and
\item $\Phi_{\upalpha_1} \circ \Phi_{\upalpha_2}^{-1} \colon \scrC_{\chamD_2} \to \scrC_{\chamD_1}$ is an atom.
\end{enumerate}
\end{defin}

\begin{remark}
The notion of atoms and their partial order can be more clearly explained using the language of the Deligne groupoid; for details see e.g.\ \cite{HaraWemyss}.
\end{remark}

The following is the key proposition, which provides the induction step.

\begin{prop} \label{HaraWemyss key prop}
Let $f \colon X \to \Spec R$ be as in Setup~\ref{geom setup},
and $x \in \scrC_{\chamC}$ an object with $\Hom(x,x[i]) = 0$ for $i < 0$.
Suppose that $x \in \llsq a,b\rrsq $ for $a < b$ and consider the set
\[ 
\nabla(x) \coloneq \{\Phi_{\upalpha} \mid \text{$\Phi_{\upalpha}$ is an atom and $\Phi_{\upalpha}(x) \in [a,b]$} \}. 
\]
\begin{enumerate}
\item\label{HaraWemyss key prop 1} 
If $\Phi_{\upalpha} \in \nabla(x)$ satisfies $\H^a(\Phi_{\upalpha}x) \neq 0$,
and $\scrS_i \hookrightarrow \H^a(\Phi_{\upalpha}x)$ is a simple submodule,
then $\Phi_i \circ \Phi_{\upalpha} \in \nabla(x)$.

\item\label{HaraWemyss key prop 2} 
If $\Phi_{\upalpha} \in \nabla(x)$ is a maximal element in this set,
then $\H^a(\Phi_{\upalpha}x) = 0$ and hence $\Phi_{\upalpha}(x) \in [a+1, b]$. Thus  $\ell(\Phi_{\upalpha}x) < \ell(x)$.
\end{enumerate}
\end{prop}

\begin{remark}
It is important to remark that $\nabla(x)$ contains a non-identity element.
Indeed, if $x \in \llsq a,b]$ for $a < b$,
there exists at least one $i$ such that $\scrS_i  \hookrightarrow \H^a(x)$.
Since the identity $\mathrm{id} \colon \scrC_{\chamC} \to \scrC_{\chamC}$ clearly belongs to $\nabla(x)$, Proposition~\ref{HaraWemyss key prop}\eqref{HaraWemyss key prop 1} shows that $s_i = s_i \circ \Id \in \nabla(x)$.
\end{remark}

Note that Proposition~\ref{HaraWemyss key prop}\eqref{HaraWemyss key prop 1} implies \eqref{HaraWemyss key prop 2}. Indeed, let $\Phi_{\upalpha_1} \coloneq \Phi_{i_1}$ be a single wall-crossing contained in $\nabla(x)$, which exits by the remark above.
If this atom $\Phi_{\upalpha_1} \in \nabla(x)$ satisfies $\H^a(\Phi_{\upalpha_1}x) \neq 0$, then there is an $i_2$ such that $\Phi_{\upalpha_2} \coloneq \Phi_{i_2} \circ \Phi_{\upalpha_1} \in \nabla(x)$.
Repeating this process gives an ascending chain
\[ 
\mathrm{id}_{\scrC_{\chamC}} < \Phi_{\upalpha_1} < \Phi_{\upalpha_2} < \hdots \in \nabla(x). 
\]
However, since $\nabla(x)$ is a finite poset, this process cannot continue forever, and hence there exits some $\upalpha_k \in \nabla(x)$
such that $\Phi_{\upalpha_k}(x) \in [a+1, b]$.

\subsection{From Bricks to Simples}

This section outlines the proof of Theorem~\ref{HaraWemyss main thm}\eqref{HaraWemyss main thm 2}.
Let $x \in \scrC_{\chamC}$ be an object such that $\Hom(x,x[i]) = 0$ for all $i < 0$ and $\Hom(x,x) = \bC$.
Thanks to \ref{HaraWemyss main thm}\eqref{HaraWemyss main thm 1},
it is enough to show the following.

\begin{prop} \label{from brick to simple}
Let $f \colon X \to \Spec R$ be as in Setup~\ref{geom setup},
and $\chamC$ a chamber of the associated hyperplane arrangement with corresponding finite-dimensional algebra $\Lambda_{\chamC, \con}$. 
Assume that an object $x \in \fmod \Lambda_{\chamC, \con}$ satisfies $\Hom(x,x) = \bC$ (i.e.~$x$ is a brick).
Then there is an atom $\Phi_{\alpha} \colon \scrC_{\chamC} \to \scrC_{\chamD}$
such that $\Phi(x) \simeq \scrO_{\Curve_i}(-1)$ for some $\Curve_i$.
\end{prop}

To prove this proposition, 
it is important to consider the realization functor of the standard heart
\[ 
F_{\chamC} \colon \Db(\fmod \Lambda_{\chamC,\con}) \to \scrC_{\chamC}, 
\]
which restricts to the identity functor $F_{\chamC}|_{\fmod \Lambda_{\chamC,\con}}\colon \fmod \Lambda_{\chamC,\con} \to \fmod \Lambda_{\chamC,\con}$ between the standard hearts.

\begin{remark}\label{rem:DbandCdifferent}
The functor $F_\chamC$ is \emph{never} an equivalence, as the categories $\scrC_{\chamC}$ and $\Db(\fmod \Lambda_{\chamC,\con})$ are very different. Nevertheless, we will end up classifying brick objects (and t-structures) in both categories: $\scrC_\chamC$ in \ref{HaraWemyss key prop} and $\Db(\fmod \Lambda_{\chamC,\con})$ as a consequence of \ref{prop:SDmain}. The answers will turn out to be the same, but there is still no good explanation of why this should be true.
\end{remark}

For an object $y \in \Db(\fmod A)$ of the derived category of an algebra $A$,
with respect to the standard heart $\fmod A$ we will recycle notation and write $y \in [a,b]$ in the obvious way, generalizing \ref{defin:abNotation}.

One of the key properties of the algebra $\Lambda_{\chamC,\con}$
is $\tau$-tilting finiteness \cite{AIR}.

\begin{defin}
Let $A$ be a finite-dimensional $\bC$-algebra.
The algebra $A$ is said to be \emph{$\tau$-tilting finite}
if the set $\{ \cS \in \smc A \mid \cS \in [-1, 0] \}$ is finite.
\end{defin}

\begin{thm} \label{thm August1}
Let $f \colon X \to \Spec R$ be as in Setup~\ref{geom setup},
and $\chamC$ a chamber of the associated hyperplane arrangement with corresponding finite-dimensional algebra $\Lambda_{\chamC, \con}$. 
Then $\Lambda_{\chamC,\con}$ is $\tau$-tilting finite.
\end{thm}

Although what is required in the proof of \ref{from brick to simple} is only the threefold case, the above theorem is true in all the geometric settings of Setup~\ref{geom setup}.

\begin{proof}
If $X = \scrZ$ is the minimal resolution of a Kleinian singularity,
then $\Lambda_{\chamC, \con}$ is a preprojective algebra of an ADE Dynkin quiver.
In this case, the $\tau$-tilting finiteness follows from \cite{Mizuno} (see also \cite{DIJ}).


If $X = \scrX$ is a threefold flopping contraction, then $\tau$-tilting finiteness of the contraction algebra $\Lambda_{\chamC, \con}$ follows from \cite{August1}.

Finally, consider the case $f\colon \scrY\to\bC^2/G$ of a partial crepant resolution of a Kleinian singularity, with associated finite-dimensional algebra $\upvarepsilon_\scrJ\Pi\upvarepsilon_\scrJ$. Using \cite[p366--367]{Pinkham} (see also \cite[\S1]{KM}), there exists a $3$-fold flopping contraction $\upphi\colon \scrX\to\Spec \scrR$ and $g\in \scrR$ such that $\Spec\scrR/g\cong\bC^2/G$, and $\upphi$ base changes to $f$. It is important that to emphasise that $g$ need not be generic. Regardless, since $\upphi$ is a flopping contraction it has an associated contraction algebra $A_{\con}$ say, and further by \cite[(3.C)]{DW1}, which does not require $g$ to be generic, necessarily $A_{\con}/g \cong \upvarepsilon_\scrJ\Pi\upvarepsilon_\scrJ$.

Via the reduction theorems \cite[4.1]{EJR} or \cite[1.5]{Kimura}, the property of being $\tau$-tilting finite descends to (and can also be detected on) central quotients.  Thus, since $A_{\con}$ is $\tau$-tilting finite, so too is the central quotient $\upvarepsilon_\scrJ\Pi\upvarepsilon_\scrJ$.
\end{proof}

Moreover, depending on the precise setting within Setup~\ref{geom setup}, the algebras $\Lambda_{\chamC, \con}$ either satisfy, or are expected to satisfy, a stronger property called silting-discrete (see \ref{ex:SDlist}\eqref{ex:SDlist 7} and \ref{rem:PreProjmaybeSD}).

The following is a key result of Asai.

\begin{thm}[\cite{Asai}] \label{thm Asai}
Let $A$ be a $\tau$-tilting finite algebra,
and $x \in \fmod A$ a semibrick module.
Then there exists an algebraic heart $\scrA \subset \Db(\fmod A)$ such that $x \in \scrA$ is a direct sum of simple objects, and further any $y \in \scrA$ satisfies $y \in [-1,0]$.
\end{thm}
We now restrict to the threefold setting of Setup~\ref{geom setup}\eqref{geom setup 3}, as this allows us to use the following theorem of August.
\begin{thm}[\cite{August}] \label{August thm2}
Let $f \colon \scrX \to \Spec \scrR$ be a threefold flopping contraction as in Setup~\ref{geom setup}\eqref{geom setup 3},
and $\scrA \subset \Db(\fmod \Lambda_{\chamC, \con})$ the heart of a bounded $t$-structure.
Then there exists
\begin{itemize}
\item a composition of mutation functors $\Phi_{\upalpha} \colon \scrC_\chamC \to \scrC_\chamD$ and
\item a standard equivalence $\Psi_{\upalpha} \colon \Db(\fmod \Lambda_{\chamC, \con}) \to \Db(\fmod \Lambda_{\chamD, \con})$
\end{itemize}
such that
\begin{enumerate}
\item[(1)] $\Psi_{\upalpha}(\scrA) = \fmod \Lambda_{\chamD, \con}$, and
\item[(2)] the following diagram of functors commutes.
\[ \begin{tikzcd}
\Db(\fmod \Lambda_{\chamC, \con}) \arrow[d, "\Psi_{\upalpha}"'] \arrow[r, "F_{\chamC}"] & \scrC_{\chamC} \arrow[d, "\Phi_{\upalpha}"] \\
\Db(\fmod \Lambda_{\chamD, \con}) \arrow[r, "F_{\chamD}"] & \scrC_{\chamD}
\end{tikzcd}
\]
\end{enumerate}
\end{thm}

Given this, the proof of \ref{from brick to simple} in the $3$-fold case proceeds as follows.  Since we can assume that $x \in \fmod \Lambda_{\chamC,\con} \subset \scrC_{\chamC}$, there exits $y \in \fmod\Lambda_{\chamC,\con} \subset \Db(\fmod \Lambda_{\chamC,\con})$ such that $F_{\chamC}(y) = x$.
Combining \ref{thm August1} and \ref{thm Asai}, there is a
heart $\scrA \subset \Db(\fmod \Lambda_{\chamC, \con})$ of a bounded $t$-structure such that $y \in \scrA$ is a simple object.
Applying \ref{August thm2} then yields a commutative diagram 
\[ \begin{tikzcd}
\Db(\fmod \Lambda_{\chamC, \con}) \arrow[d, "\Psi_{\upalpha}"'] \arrow[r, "F_{\chamC}"] & \scrC_{\chamC} \arrow[d, "\Phi_{\upalpha}"] \\
\Db(\fmod \Lambda_{\chamD, \con}) \arrow[r, "F_{\chamD}"] & \scrC_{\chamD}
\end{tikzcd}
\]
such that $\Psi_{\upalpha}(\scrA) = \fmod \Lambda_{\chamD,\con}$.
Since $\Psi_{\upalpha}(y) \in \fmod \Lambda_{\chamD,\con}$ is a simple object, the commutativity of the diagram and \ref{rem simples of std heart} shows that 
\[
\Phi_{\upalpha}(x) \cong \Phi_{\upalpha}F_{\chamC}(y) \cong F_{\chamD}\Psi_{\upalpha}(y) \cong  \scrS_j
\] 
for some simple $\Lambda_\chamD$-module $\scrS_j$. As already remarked in \ref{rem simples of std heart}, $\scrS_j$ corresponds to $\scrO_{\Curve_j}(-1)$ across the derived equivalence with $X_{\chamD}$.
This finishes the proof of the $3$-fold setting within Setup~\ref{geom setup}\eqref{geom setup 3}.

The remaining cases of Setup~\ref{geom setup} in dimension two are proved by reduction to the $3$-fold situation through deformations; see \cite[3.6]{HaraWemyss}. This step is mildly subtle, since \ref{August thm2} is a result in dimension three, which is false  in dimension two.

\subsection{Classification of $t$-structures}

Our classification technique of spherical objects extends to give the following classification of all bounded $t$-structures of the null category $\scrC$.

\begin{thm}
Let $f \colon X \to \Spec R$ be as in Setup~\ref{geom setup},
$\chamC$ a chamber of the associated hyperplane arrangement,
and $\scrA \subset \scrC_{\chamC}$ the heart of a bounded $t$-structure.
Then there is a composition of mutation functors
\[ 
\Phi_{\upalpha} \colon \scrC_{\chamC} \to \scrC_{\chamD} 
\]
such that $\Phi_{\upalpha}(\scrA) = \fmod \Lambda_{\chamD,\con}$ is the standard heart.
In particular, all bounded hearts of $\scrC_{\chamC}$ are algebraic.
\end{thm}

From this, it follows that for any of the settings in Setup~\ref{geom setup}, the associated null category $\scrC$ is \emph{$t$-discrete} in the sense of \cite{AMY}.

\section{`Finite' Geometric Setting \`{a} la Bapat--Deopurkar--Licata}

The paper \cite{BDL} of Bapat, Duopurkar, and Licata classifies spherical objects in the category $\scrC$ in Setup~\ref{geom setup}\eqref{geom setup 1}, namely minimal resolutions of Kleinian singularities,  using the notion of Bridgeland stability conditions \cite{Bstab}.

Let $\scrT$ be a $2$-Calabi-Yau category with finite dimensional $\Hom$-spaces, equipped with a Bridgeland stability condition $\sigma$.
Every object $x \in \scrT$ has a Harder–Narasimhan (HN) filtration
\[
\begin{tikzcd}
0 = x_0 \arrow[rr] & & x_1 \arrow[dl] \arrow[r] & \hdots \arrow[r] & x_{n-1} \arrow[dl] \arrow[rr] & & x_n = x, \arrow[dl] \\
& z_1 \arrow[lu, dashrightarrow] & & \hdots \arrow[lu, dashrightarrow] & & z_n \arrow[lu, dashrightarrow] &
\end{tikzcd}
\]
where each $z_i$ is $\sigma$-semistable, and their phases satisfy
\[ \upphi(z_1) > \hdots > \upphi(z_n). \]
The top and bottom phases are denoted as $\upphi_+(x) \coloneq \upphi(z_1)$
and $\upphi_-(x) \coloneq \upphi(z_n)$.
The $\sigma$-semistable objects $z_i$ are filtered by $\sigma$-stable objects in the same phases by Jordan-H\"{o}lder filtrations.
In this strategy, the key proposition is the following.

\begin{prop}[{\cite[2.3]{Huy}, \cite[4.1]{BDL}}]
Let $\scrT$ and $\sigma$ be as above.
\begin{enumerate}
\item[(1)] If $x \in \scrT$ is a spherical object, then any stable factor of $x$ is again spherical.
\end{enumerate}
Assume that $\scrT$ is equivalent to the null subcategory $\scrC \subset \Db(\coh \scrZ)$ of the minimal resolution $f \colon \scrZ \to \bC^2/G$ for some finite subgroup $G \leq \SL(2,\bC)$.
Suppose in addition that the stability $\sigma$ is generic in the sense of \cite[4.1]{BDL}.
Then
\begin{enumerate}
\item[(2)] any semistable factor $z_i$ of a spherical object $x \in \scrC$
is isomorphic to $y_i^{\oplus m_i}$ for a $\sigma$-stable spherical object $y_i$ and $m_i > 0$.
\end{enumerate}
\end{prop}

The key idea in \cite{BDL} is to reduce the spread 
\[ \ell(x) \coloneq \upphi_+(x) - \upphi_-(x) \]
of phases by applying spherical twists.

From now on, consider the null subcategory $\scrC$ of the minimal resolution $f \colon \scrZ \to \bC^2/G$ for some finite subgroup $G\leq \SL(2,\bC)$.
Consider the standard heart $\fmod \Lambda_{\con} \subset \scrC$,
and let $\sigma = (Z, \fmod \Lambda_{\con})$ be a generic Bridgeland stability condition on $\scrC$.
In this setup, if the spread $\ell(x)$ of phases attains the minimal value, i.e.~if $\ell(x) = 0$, then $x$ is a $\sigma$-stable object.
Such an object can be classified as in \ref{from brick to simple}.
Namely, if $x$ is as above, then there are simple modules $s_0, s_1, \hdots, s_t$ in $\fmod \Lambda_{\con}$ such that
\begin{align} \label{stable to simple}
x \simeq T_{s_t} \circ \hdots T_{s_1}(s_0).
\end{align}
Note that each $s_i$ is isomorphic to $\scrO_{\Curve_{j_i}}(-1)$ for some 
irreducible component $\Curve_{j_i}$.

The step that makes $\ell(x)$ smaller is  the following.

\begin{thm} \label{BDL improve result}
In the setup as above,
let $x \in \scrC$ be a spherical object with $\ell(x) > 0$, 
with HN filtration given by
\[
\begin{tikzcd}
0 = x_0 \arrow[rr] & & x_1 \arrow[dl] \arrow[r] & \hdots \arrow[r] & x_{n-1} \arrow[dl] \arrow[rr] & & x_n = x, \arrow[dl] \\
& z_1 \arrow[lu, dashrightarrow] & & \hdots \arrow[lu, dashrightarrow] & & z_n \arrow[lu, dashrightarrow] &
\end{tikzcd}
\]
If $y_1$ and $y_n$ are $\sigma$-stable spherical objects such that
$z_i \simeq y_i^{\oplus m_i}$, then
\begin{enumerate}
\item $\ell(T_{y_1}(x))<\ell(x)$.
\item $\ell(T_{y_n}^{-1}(x))<\ell(x) $.
\end{enumerate}
\end{thm}

For the null subcategory $\scrC$ and any positive real number $\varepsilon \in \bR_{>0}$,
the set
\[ \{\ell(x') \in \bR \mid \text{$x' \in \scrC$ is spherical and $\ell(x') \leq \varepsilon$} \} \subset \bR \]
is a finite set. (Note that this fact does not hold for the category $\scrD$ of complexes with compact support!) Therefore, applying the theorem repeatedly and then using (\ref{stable to simple}) gives the main result.

\begin{cor}[{\cite[1.1]{BDL}}] \label{BDL classification}
For any spherical object $x \in \scrC$, there is
\[
\Psi \in B' = \langle T_{s} \mid \text{$s \in \fmod \Lambda_{\con}$ is a $\sigma$-stable spherical object} \rangle \subset \Auteq \scrC 
\]
such that $\Psi(x) \simeq \scrO_{\Curve_j}(-1)$
for some irreducible curve $\Curve_j$.
\end{cor}

Note that (\ref{stable to simple}) also implies that
\[ B' = \left\langle T_{s_j} \mid \text{$s_j \coloneq \scrO_{\Curve_j}(-1)$ for an irreducible curve $\Curve_j \subset \Curve_{\mathrm{red}}$} \right\rangle. \]

\begin{remark}
Note that \ref{BDL improve result} requires $x \in \scrC$ to be spherical,
whilst \ref{HaraWemyss main thm}\eqref{HaraWemyss main thm 2} holds for weaker objects $x \in \scrC$ such that $\Hom(x,x[i]) = 0$ for all $i < 0$ and $\Hom(x,x) = \bC$ (i.e.~$x$ is a brick). It is remarkable that the two results have different assumptions on $x \in \scrC$, but they have the same classification.
\end{remark}

\section{Silting-discrete Setting}\label{sec:SiltingDiscrete}

The main technique in \S\ref{subsec:1.4CLass} also applies to the derived category of silting-discrete algebras in representation theory.
In this section, we review this work as well as the works by Aihara--Mizuno \cite{AM},
Pauksztello--Saor\'{i}n--Zvonareva \cite{David},
and Adachi--Mizuno--Yang \cite{AMY}.

\subsection{Silting Complexes and the K\"{o}ning--Yang Bijection}

As is clear from the definition in \ref{def:silting}, the notion of silting complexes is a generalization of tilting complexes, and hence of projective generators.
Note that, in contrast to tilting complexes, a silting complex is allowed to have negative self-extensions.

A key property of silting complexes is that they admit a natural partial order.

\begin{defin}\cite{KY}
Let $A$ be a finite dimensional algebra over $\bC$.
\begin{enumerate}
\item For silting complexes $P, Q \in \Kb(\proj A)$,
the notation $P \geq Q$ indicates that $\Hom(P,Q[i]) = 0$ for all $i > 0$.
\item For smcs $x, y \in \Db(\fmod A)$,
the notation $x \geq y$ indicates that $\Hom(y,x[i]) = 0$ for all $i < 0$.
\end{enumerate}
\end{defin}
Write $\silt A$ for the set of silting complexes in $\Kb(\proj A)$, up to isomorphism, and $\smc A$ for the set of simple minded collections in $\Db(\fmod A)$.  The above relations~$\geq$ define partial orders on both $\silt A$ and $\smc A$.


\begin{thm}[The K\"onig-Yang Bijection \cite{KY}] \label{KY bij}
Let $A$ be a finite dimensional algebra over $\bC$.
Then there are order-preserving bijections between the following sets.
\begin{enumerate}
\item[(1)] $\silt A$
\item[(2)] $\smc A$
\item[(3)] The set of bounded $t$-structures on $\Db(\fmod A)$ with algebraic hearts.
\end{enumerate}
Let $S = \{x_1, \hdots, x_n\}$ be an smc, and $P_{S}$ and $\scrA_{S}$ the corresponding silting complex and algebraic heart.
Then $\scrA_{S} \simeq \fmod \End(P_{S})$ and $x_1, \hdots, x_n$ gives the full collection of simple $\End(P_{S})$-modules.
\end{thm}

The above theorem roughly says that silting complexes, smcs, and bounded algebraic hearts classify each other.
Note that the K\"onig-Yang bijection does not give any control for general presilting complexes, semibricks, and bounded $t$-structures.

\subsection{Silting-discrete Classifications}
The derived category $\Db(\fmod A)$ where $A$ is a silting-discrete algebra is another class of categories in which the theory established by Iyama et.al.~works dreamily well (see, for example, \ref{silsing discrete classification} below). Such algebras are defined as follows.

\begin{defin}[\cite{AM}] \label{silting-discrete definition} 
A finite dimensional $\bC$-algebra $A$ is called \emph{silting-discrete}
if one of (and hence both of) the following two equivalent conditions are satisfied.
\begin{enumerate}
\item[(1)] For any $P \in \silt A$, the set $\{Q \in \silt A \mid P \geq Q \geq P[1] \}$ is finite.
\item[(2)] For any $x \in \smc A$, the set $\{y \in \smc A \mid x \geq y \geq x[1] \}$ is finite.
\end{enumerate}
\end{defin}

\begin{example}[{c.f.~\cite{AM, August1, AD}}]\label{ex:SDlist}
The following algebras are silting-discrete.
\begin{enumerate}
\item\label{ex:SDlist 1} A path algebra of Dynkin type.
\item\label{ex:SDlist 2} A local algebra.
\item\label{ex:SDlist 3} A representation-finite symmetric algebra.
\item\label{ex:SDlist 4} A derived discrete algebra of finite global dimension.
\item\label{ex:SDlist 5} A Brauer graph algebra whose Brauer graph contains at most one cycle of odd length and no cycle of even length.
\item\label{ex:SDlist 6} The preprojective algebra of Dynkin type $D_{2n}$, $E_7$, or $E_8$.
\item\label{ex:SDlist 7} The contraction algebra $\Lambda_{\con}$ for any $3$-fold flopping contraction $f \colon \scrX \to \Spec \scrR$ as in Setup~\ref{geom setup}\eqref{geom setup 3}. By \cite{Hao}, these include the following large class of examples.  Consider the following quiver $Q$, where $n\geq 2$.
\[
\begin{array}{c}
\begin{tikzpicture}[scale=1.25,bend angle=15, looseness=1]
\node (a) at (-1,0) [pvertex] {$1$};
\node (b) at (0,0) [pvertex] {$2$};
\node at (0,0.6) [pvertex] {};
\node (c) at (1,0) [pvertex] {$\phantom{2}$};
\node at (1,0) {$\scriptstyle \hdots$};
\node (n) at (2,0) [pvertex] {$n$};
\draw[->,bend left] (b) to node[below] {$\scriptstyle b_2$}(a);
\draw[->,bend left] (a) to node[above] {$\scriptstyle a_2$} (b);
\draw[->,bend left] (c) to node[below] {$\scriptstyle b_4$} (b);
\draw[->,bend left] (b) to node[above] {$\scriptstyle a_4$}  (c);
\draw[->,bend left] (n) to node[below] {$\scriptstyle b_{2n-2}$} (c);
\draw[->,bend left] (c) to node[above] {$\scriptstyle a_{2n-2}$}  (n);
\draw[<-]  (a) edge [in=-120,out=-65,loop,looseness=7] node[below] {$\scriptstyle \mathtt{x}_1$} (a);
\draw[<-]  (b) edge [in=-120,out=-65,loop,looseness=7] node[below] {$\scriptstyle \mathtt{x}_3$} (b);
\draw[<-]  (n) edge [in=-120,out=-65,loop,looseness=7] node[below] {$\scriptstyle \mathtt{x}_{2n-1}$} (n);
\end{tikzpicture}
\end{array}
\]
Define $\mathtt{x}_{2i}=a_{2i}b_{2i}$, and set $\mathtt{x}'_{2i+1}=\mathtt{x}_{2i+1}$ and $\mathtt{x}'_{2i}=b_{2i}a_{2i}$, then consider
\[
W=\sum_{i=1}^{2n-2}\mathtt{x}^\prime_{i}\mathtt{x}^{\phantom\prime}_{i+1}+\sum_{i=1}^{2n-1}\sum_{j=2}^\infty \upkappa_{ij}\mathtt{x}_i^j
\]
for some choice of scalars $\upkappa_{ij}$.  For any such choice for which the (completed) Jacobi algebra $\scrJ\mathrm{ac}(W)$ is finite dimensional, the algebra $\scrJ\mathrm{ac}(W)$ is the contraction algebra of a $cA_n$ smooth flopping contraction.
\end{enumerate}
\end{example}

\begin{remark}\label{rem:PreProjmaybeSD}
It is frustrating that whilst in Example~\ref{ex:SDlist}\eqref{ex:SDlist 6} the preprojective algebra of Types $D_{2n}$, $E_7$, and $E_8$ are known to be silting-discrete, Type $A_n$ in general is still open.  The cases $A_n$ with  $n = 1,2$ are known to be silting-discrete, and it seems expected that this holds in general.

If all ADE preprojective algebras are indeed silting-discrete, then it follows immediately from \cite{AH2} that the induced contracted preprojective algebras $\upvarepsilon_{\scrJ} \Uppi \upvarepsilon_{\scrJ}$ are also silting-discrete.
\end{remark}

One beautiful aspect of silting-discrete algebras is that presilting complexes, semibricks, and $t$-structures have the following properties.

\begin{thm} \label{silsing discrete classification}
Let $A$ be a silting-discrete $\bC$-algebra.
Then the following hold.
\begin{enumerate}
\item\label{silsing discrete classification 1} \cite{AM} Let $P$ be a presilting complex in $\Kb(\proj A)$.
Then there exists $P' \in \Kb(\proj A)$ such that $P \oplus P'$ is silting.
\item\label{silsing discrete classification 2} \cite{HaraWemyss} Let $x$ be a semibrick complex in $\Db(\fmod A)$.
Then there exists $x' \in \Db(\fmod A)$ such that $x \oplus x'$ is an smc. 
\item\label{silsing discrete classification 3} \cite{David, AMY} Let $\scrA$ be the heart of a bounded $t$-structure on $\Db(\fmod A)$.
Then $\scrA$ is algebraic, hence $\scrA \simeq \fmod B$ for some finite dimensional $\bC$-algebra $B$.
\end{enumerate}
\end{thm}

As remarked in \ref{rem:Dual}, silting and smcs are in some sense dual, but to the best of the authors' knowledge, there are no implications between the statements, or proofs, of (1) and (2) above.  


In the following, let us outline the proof of \ref{silsing discrete classification}\eqref{silsing discrete classification 2}, which is roughly divided into the following two steps.
\begin{enumerate}
\item[(a)] Find an algebraic heart $\scrA \subset \Db(\fmod A)$ such that $x \in \scrA$.
\item[(b)] Find a complex $x'$ which is two-term with respect to the heart $\scrA$ such that $x \oplus x'$ is an smc.
\end{enumerate}
Step (b) is quite easy, as it uses the following lemma which is a direct consequence of the definition.
\begin{lemma}
Let $A$ be silting-discrete, $\scrA \subset \Db(\fmod A)$ an algebraic heart,
and $B$ an algebra such that $\scrA \simeq \fmod B$.
Then $B$ is $\tau$-tilting finite. In particular, by \ref{thm Asai}, any semibrick $x \in \scrA$ can be completed into an smc that is two-term with respect to $\scrA$.
\end{lemma}

All the work goes into Step (a), and perhaps a little surprisingly, this step works for more general objects than semibricks, namely those objects $x \in \Db(\fmod A)$ such that $\Hom(x,x[i]) = 0$ for all $i < 0$.

\begin{remark}\label{rem:notenoughAut}
The proof of step (a) is very similar to the proof of \ref{HaraWemyss main thm}\eqref{HaraWemyss main thm 2}, but the main difference is that for a general silting-discrete algebra there are few autoequivalences. Consequently the induction step uses simple-minded mutations to reduce the homological spread, in lieu of the mutation functors before.
\end{remark}

\begin{defin}
Let 
 $\cS = \{x_1, x_2, \hdots, x_n \}$ an smc in $\Db(\fmod A)$.
For a chosen $x_i \in \cS$, define a new collection
$\cS' = \{x_1', \dots, x_n'\}$ as follows.
\begin{enumerate}
\item[(1)] $x_i' \coloneq x_i[1]$.
\item[(2)] If $j \neq i$, then let $x_j[-1] \to x_{ij}$ be the minimal left approximation in the extension closure of $x_i$,
and put $x_j' \coloneq \mathrm{Cone}(x_j[-1] \to x_{ij})$.
\end{enumerate}
The collection $\cS'$ is again an smc \cite[7.6]{KY},
called the \emph{(left) mutation} of $\cS$ at $x_i$,
and is denoted $\upnu_i\cS$.
\end{defin}

Notions of mutations for silting complexes and $t$-structures also exist \cite{AI}, and all are compatible under the K\"onig-Yang bijection \ref{KY bij}.
In order to distinguish them, mutations of smcs are sometimes called simple-minded mutations.

The first step for proving \ref{silsing discrete classification}\eqref{silsing discrete classification 2}
is to show that a given semibrick is contained in an algebraic heart.
For this, let us introduce the following notation, which is analogous to \ref{defin:abNotation}.

\begin{defin}
For each smc $\cS$ in $\Db(\fmod A)$, where $A$ is a finite dimensional algebra, let $\scrA_{\cS} \subset \Db(\fmod A)$ denote the corresponding algebraic heart under \ref{KY bij}, and let $\H_{\cS}^* \colon \Db(\fmod A) \to \scrA_{\cS}$ be the associated cohomology functor. Then, for $x \in \Db(\fmod A)$ and integers $a \leq b$,
\begin{enumerate}
\item the notation $x \in [a, b]_{\cS}$ indicates that
$\H_{\cS}^i(x) = 0$ for all $i < a$ and $i>b$.

\item the notation $x \in \llsq a, b \rrsq _{\cS}$ indicates that $x \in [a,b]$ and
$\H_{\cS}^i(x) \neq 0$ for $i = a, b$.
\end{enumerate}
If $x \in \llsq a, b \rrsq _{\cS}$, put
\[ \ell(x, \cS) \coloneq b - a. \]
\end{defin}

Note that $\ell(x, \cS) = 0$ if and only if $x \in \scrA_{\cS}[k]$ for some $k \in \bZ$. Therefore, in order to establish that Step (a) is true,
it suffices to show that there always exists another smc $\cS'$ such that $\ell(x, \cS')<\ell(x, \cS)$. This is achieved in the following.

\begin{prop}\label{prop:SDmain}
Let $A$ be a silting-discrete algebra, 
and $x \in \Db(\fmod A)$ an object such that $\Hom(x,x[i]) = 0$ for all $i < 0$.
Let $\cS$ be an smc and assume that $x \in \llsq a,b\rrsq $ for some $a < b \in \bZ$.
Consider the set
\[ \nabla(x, \cS) \coloneq \{ \cS' \in \smc A \mid \text{$\cS \geq \cS' \geq \cS[1]$ and $x \in [a,b]_{\cS}$} \} \subset \smc A,  \]
which is finite since $A$ is silting-discrete.

\begin{enumerate}
\item[(1)] Assume that $\cS' = \{y_1', \hdots, y_n' \} \in \nabla(x, \cS)$ satisfies $\H_{\cS'}^a(x) \neq 0$. 
Then for a simple subobject $y'_i \hookrightarrow \H_{\cS'}^a(x)$ in the heart $\scrA_{\cS'}$,
the mutation $\upnu_i \cS'$ satisfies 
$\upnu_i \cS' \in \nabla(x, \cS)$ and $\cS' > \upnu_i \cS'$.

\item[(2)] If $\cS'$ is a minimal element of $\nabla(x, \cS)$,
then $H_{\cS'}^a(x) = 0$ and hence $x \in [a+1, b]$. Thus $\ell(x, \cS') < \ell(x,\cS)$.
\end{enumerate}
\end{prop}

Applying the above repeatedly shows that following, which establishes Step (a).

\begin{cor}
Let $A$ be a silting-discrete algebra, 
and $x \in \Db(\fmod A)$ an object.
Then the following are equivalent.
\begin{enumerate}
\item[(1)] $\Hom(x,x[i]) = 0$ for all $i < 0$.
\item[(2)] Then there exists an algebraic heart $\scrA \subset \Db(\fmod A)$ (and hence an smc $\cS$ with $\scrA = \scrA_{\cS}$) such that $x \in \scrA$.
\end{enumerate}
\end{cor}

\subsection{Threefold Contraction Algebras}
The above is for general silting-discrete algebras; for contraction algebras the autoequivalence group is larger, so more can be said.  Reinterpreting \ref{silsing discrete classification}\eqref{silsing discrete classification 2} gives the following, and the most remarkable aspect is that this is basically identical to \ref{HaraWemyss main thm}, with the caveat that $\Db(\fmod\Lambda_{\con})$ and $\scrC$ are very different categories (see \ref{rem:DbandCdifferent}).
\begin{thm} \label{HWContAlgMain}
Consider the contraction algebra $\Lambda_{\chamC, \con}$ associated to a $3$-fold flopping contraction $f_\chamC \colon \scrX_\chamC \to \Spec \scrR$.
Then the following statements hold.
\begin{enumerate}
\item\label{HWContAlgMain 1} An $x \in \Db(\fmod\Lambda_{\chamC, \con})$ satisfies $\Hom(x,x[i]) = 0$ for all $i < 0$ if and only if there is a composition of standard equivalences $\Psi_{\upalpha}\colon\Db(\fmod\Lambda_{\chamC, \con})\to\Db(\fmod\Lambda_{\chamD, \con})$, for some chamber $\chamD$, such that $\Psi_{\upalpha}(x)$ is contained in the standard heart $\fmod \Lambda_{\chamD, \con}$.
\item\label{HWContAlgMain 2} An $x \in \Db(\fmod\Lambda_{\chamC, \con})$ satisfies $\Hom(x,x[i]) = 0$ for all $i < 0$ and further $\Hom(x,x) = \bC$ if and only if there is a composition of standard equivalences $\Psi_{\upalpha}\colon\Db(\fmod\Lambda_{\chamC, \con})\to\Db(\fmod\Lambda_{\chamD, \con})$, for some chamber $\chamD$, such that $\Psi_{\upalpha}(x)\cong \scrS_j$ for some simple module $\scrS_j$.
\end{enumerate}
\end{thm}

\section{Symplectic Setting \`{a} la Smith--Wemyss}

This section summarises some results in symplectic geometry obtained via the algebraic geometry setup above, discovered in \cite{SW}, mainly since in this restricted class of examples the algorithm that classifies spherical objects is significantly more direct. 

As motivation, it is well-known that the mirror to the minimal resolution of a Type $A$ Kleinian singularity can be obtained by plumbing a chain of $S^2$'s along points. In the case of $A_2$, the cartoon picture is simply
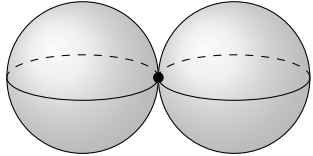
\[
\begin{tikzpicture}[scale=0.5]
  \shade[ball color = gray!40, opacity = 0.4] (0,0) circle (2cm);
  \draw (0,0) circle (2cm);
  \draw (-2,0) arc (180:360:2 and 0.6);
  \draw[dashed] (2,0) arc (0:180:2 and 0.6);
  
    \shade[ball color = gray!40, opacity = 0.4] (4,0) circle (2cm);
  \draw (4,0) circle (2cm);
  \draw (2,0) arc (180:360:2 and 0.6);
  \draw[dashed] (6,0) arc (0:180:2 and 0.6);
  \fill[fill=black] (2,0) circle (4pt);
\end{tikzpicture}
\]
To mirror 2-curve flops is significantly harder. 

We now restrict the setting to the very specific case of a flopping contraction with two curves, in which the normal bundle of both curves is $(-1,-1)$.  The latter assumption ensures that the structure sheaves of the curves are spherical objects.  Such flopping contractions are classified by finite dimensional Jacobi algebras on the quiver
\[
Q=\begin{array}{c}
\begin{tikzpicture}[scale=0.8]
\node (A) at (0,0) [B] {};
\node (B) at (2,0) [B] {};
\draw[->, bend left] (A) to node[above]{$\scriptstyle a$} (B);
\draw[->, bend left] (B) to node[below]{$\scriptstyle b$} (A);
\end{tikzpicture}
\end{array}
\]
for which \cite[3.6]{DWZ} asserts that, after change in coordinates, the potential has the form $\mathbb{W}_k= (ab)^k$ for some $k\geq 1$. Thus necessarily the flopping contraction is Type $cA_2$, and up to analytic change in coordinates (over $\bC$) the base is given by the single equation $uv-xy(x^k+y)$, where $k\geq 1$.  We will write this flopping contraction as $f_k\colon \scrX_k\to\Spec \scrR_k$, to emphasise the dependence on $k$.

\medskip
To mirror this requires two spherical objects (so, spheres!) to be plumbed in some way. This is achieved via the \emph{double bubble plumbing}, namely the Stein manifold obtained by plumbing the cotangent bundles of two 3-spheres $Q_0$ and $Q_1$ along an unknotted circle $Z\subset Q_i$. Strictly speaking, this requires an identification 
\[
\upeta\colon \upnu_{Z/Q_0} \stackrel{\sim}{\longrightarrow} \upnu_{Z/Q_1},
\]
and so in fact we obtain a family of such plumbings $W_k$, indexed over $k\in\mathbb{Z}$.  Since $W_k\cong W_{-k}$, the family depends only on $k\geq 0$, and since the case $k=0$ behaves very differently, henceforth we will assume that $k\geq 1$.

\medskip
Now there is a nullcategory $\scrC_k$ associated to the flopping contraction $f_k$, and a Fukaya category $\scrF\mathrm{uk}(W_k)$ associated to $W_k$.  Working over the field $\bC$ does not line up these categories correctly \cite[1.1(1)]{SW}, but over a field $\mathbb{K}$ of characteristic $p=k$, with appropriate tweaks including setting up the birational geometry over $\mathbb{K}$, the categories  $\scrC_k$ and $\scrF\mathrm{uk}(W_k)$ turn out to be equivalent \cite[1.1(2)]{SW}.

The following quirk in the representation theory of $\Lambda_{\con}=\scrJ\mathrm{ac}(Q,\mathbb{W}_k)$ allows us to give a much more direct route towards classification.
\begin{prop}\label{prop;SWsimplify}
Consider $f_k\colon \scrX_k\to\Spec \scrR_k$, with associated $\Lambda_{\con}$.
 If $M,N\in\fmod\Lambda_{\con}$ with $\Hom(M,N)=0$, then either $M$ or $N$ is filtered by a unique simple.  
\end{prop} 
Given an object $x$ with $\Hom_{\scrC_k}(x,x[i])=0$ for all $i<0$, say with homological spread lying in $\llsq a,b\rrsq $ with $a<b$, a spectral sequence shows that $\Hom(\mathrm{H}^bx,\mathrm{H}^ax)=0$. By \ref{prop;SWsimplify} there are four options, and below each option we assign the following functor.
\[
\begin{array}{rcccc}
&\mathrm{H}^ax&\mathrm{H}^ax&\mathrm{H}^bx&\mathrm{H}^bx\\
\mbox{filtered by:}&\scrS_1&\scrS_2&\scrS_1&\scrS_2\\[2mm]
\mbox{functor:}&\Phi_1&\Phi_2&\Phi_1^{-1}&\Phi_2^{-1}\\
\end{array}
\]
The following asserts that the above assignment always improves the homological spread, and so gives a very direct induction step to the classification of objects $x$ satisfying $\Hom_{\scrC_k}(x,x[i])=0$ for all $i<0$.  This should be contrasted with \S\ref{subsec:1.4CLass} where in general testing against compositions of positive mutation functors is required.  Of course, the point is that whilst testing against atoms in general is required, in some situations (like here) more direct approaches can work.
\begin{thm}
Suppose that $ x\in\scrC_k$  satisfies $\Hom_{\scrC_k}( x, x[i])=0$ for all $i<0$, and $x\in\llsq a,b\rrsq $ with $\ell(x)=b-a>0$. Then the following statements hold. 
\begin{enumerate}
\item Either $\mathrm{H}^ax$ or $\mathrm{H}^bx$ is filtered by a unique simple module $\scrS_i$, with $i\in\{1,2\}$.
\item When $\mathrm{H}^a(\Psi a)$ is filtered only by $\scrS_i$, applying $\Phi_i$ decreases the length $\ell$.
\item When $\mathrm{H}^b(\Psi a)$  is filtered only by $\scrS_i$, applying $\Phi_i^{-1}$ decreases the length $\ell$.
\end{enumerate}
\end{thm}

The above result provides part of the induction step which leads to the classification of various classes of objects, exactly as in \ref{HaraWemyss main thm}. There are various symplectic corollaries, including the following.

\begin{cor}
Let $p=1$ or $p>2$ be prime. If $L\subset W_p$ is a closed exact Lagrangian submanifold with vanishing Maslov class, then $\pm [L] \in \{(1,0), (0,1), (1,\pm 1)\} \subset \bZ\oplus\bZ$.
\end{cor}

\section{`Affine' Geometric Setting (surfaces) \`{a} la Ishii--Uehara}

Here, consider the minimal resolution $f \colon \scrZ \to \bC^2/G$,
where $G \subset \SL(2,\bC)$ is a cyclic subgroup; in other words, $\bC^2/G$ is a Type $A$ Kleinian singularity. As before, let $\Curve = f^{-1}(0)$ be the scheme-theoretic fiber above the unique singular point, and now consider the `affine' category 
\[
\scrD=\{ a\in\Db(\coh \scrZ)\mid \Supp a\subseteq \Curve\}.
\]
Write $\cH_{\coh} \colon \scrD \to \coh_\Curve(Z)$ for the cohomology functor with respect to the heart $\coh_\Curve(Z)\subset\scrD$. The calligraphic font on $\cH_{\coh}$ reminds us that the output is now a sheaf, not a module.

This setting is significantly harder than the categories $\scrC$ considered before. One of the main observations in \cite{IU} is the following, which asserts that coherent cohomology of spherical objects is well behaved.

\begin{prop}[{\cite[4.10, 6.1]{IU}}]\label{prop:IUtech}
Let $x \in \scrD$ be a spherical object.
Then the following statements hold.
\begin{enumerate}
\item\label{prop:IUtech 1} $\bigoplus_{i \in \bZ} \cH_{\coh}^i(x)$ is a rigid $\scrO_\Curve$-module, pure of dimension one.
\item\label{prop:IUtech 2} Furthermore, each indecomposable summand of $\cH_{\coh}^i(x)$
is a line bundle on a connected pure one-dimensional subscheme of $\Curve$.
\end{enumerate}
\end{prop}

Note that (1) applies to an arbitrarily finite subgroup $G \leq \SL(2, \bC)$, however (2) heavily relies on the assumption that the singularity of $\bC^2/G$ is of Type $A$.

Let $\Curve_1, \Curve_2, \hdots, \Curve_n$ be the irreducible components of $\Curve_{\mathrm{red}}$ such that $\Curve_{\mathrm{red}} = \bigcup \Curve_i$,
and $\eta_i \in \Curve_i$ the generic point.
For a spherical object $x \in \scrD$, the key invariant here is 
\[ 
\ell(x) \coloneq \sum_{i, p} \mathrm{length}_{\scrO_{\Curve_i, \eta_i}} \cH^p_{\coh}(x)_{\eta_i}. 
\]
Notice that if $\ell(x) = 1$, then $x \simeq \scrO_{\Curve_i}(a)[k]$ for some $i, a, k$.
The following is then the induction step, which implies the classification of spherical objects.

\begin{thm}[\cite{IU}] \label{IU main}
For any spherical object $x \in \scrD$ with $\ell(x) > 1$,
there exists an element
\[ \Phi \in B \coloneq \langle T_{\scrO_{\Curve_i}(-1)}, T_{\omega_\Curve} \mid \Curve_i \subset \Curve \rangle \subset \Auteq \scrD \]
such that $\ell(\Phi x)<\ell(x)$. Consequently, there is $\Psi \in B$ such that $\Psi(x) \simeq \scrO_{\Curve_i}(a)[k]$ for some $i, a, k$.
\end{thm}

It seems difficult to generalize the method in \cite{IU} to arbitrary Kleinian singularities, partly since \ref{prop:IUtech}\eqref{prop:IUtech 2} so heavily relies on the Type $A$ assumption (c.f.~\cite{Kawatani}).
Furthermore, within \cite{IU} is a key spectral sequence which heavily relies on two facts: the first is that the dimension is two, and the second that the objects $y$ satisfy $\Hom(y,y[1])=0$. As such, the spectral sequence cannot be applied to classify more general objects $y$ which only satisfy $\Hom(y,y[i]) = 0$ for all $i < 0$, nor can it be applied to the 3-fold flops setting.

At the time of writing, there are very few works that directly apply to the `affine' category $\scrD$ of compactly support complexes. In addition to \cite{IU} are the papers \cite{KeatingSmith, Parth1, Parth2}, which are also surveyed below. 

Note that the induction step in \cite{IU} directly descends to the induction in the null subcategory $\scrC$. In particular, the proof of \ref{IU main} contains the proof of the following, which was first explicitly stated in \cite{IUU}.

\begin{cor}[\cite{IUU}] \label{IUU finite classification}
For any spherical object $x \in \scrC$,
there exists an element
\[ 
\Phi \in B' \coloneq \langle T_{\scrO_{\Curve_i}(-1)} \mid \Curve_i \subset \Curve \rangle \subset \Auteq \scrC 
\]
such that $\ell(\Phi x)<\ell(x)$. Consequently, there is $\Psi \in B'$ such that $\Psi(x) \simeq \scrO_{\Curve_i}(-1)$ for some $i$.
\end{cor}

As in \ref{BDL classification}, this result \ref{IUU finite classification} is now superseded by \ref{HaraWemyss main thm}, but again the proofs are very different.

\begin{remark}\label{rem:Surj}
Let $f \colon \scrZ \to \bC^2/G$ be the minimal resolution of a Type A Kleinian singularity.
As a consequence of the classification result \ref{IU main} of spherical objects, it follows that
\[ \langle T_{\scrS} \mid \text{$\scrS \in \scrD$  is spherical} \rangle = B \coloneq  \langle T_{\scrO_{\Curve_i}(-1)}, T_{\omega_\Curve} \mid \Curve_i \subset \Curve \rangle \subset \Auteq \scrD. \]
As a corollary, it is shown in \cite[1.3]{IU} (see also \cite[Appendix~A]{IUU}) that
\[ 
\Auteq \scrD = (\langle B, \Pic X \rangle \rtimes \Aut X) \times \bZ. 
\]
\end{remark}

\begin{remark}
The classification of spherical objects and the resulting description of $\Auteq\scrD$ in \ref{rem:Surj} is related to the question of surjectivity of the homomorphism from the (affine) braid group to the subgroup of the autoequivalence group generated by spherical twists. The question of injectivity is also interesting, and difficult in general.  For the minimal resolution of a Type A Kleinian singularity, injectivity has been proved in \cite[Cor.~37]{IUU},
and hence the group $B$ above is isomorphic to the affine braid group of Type A.
It is remarkable that the proof in \cite{IUU} is achieved by reduction to characteristic two, a technique which is not possible in the $3$-fold setting, where the defining equations do not in general have integer coefficients.  We remark that injectivity in the `finite' settings of the null category $\scrC$ is easier, and is contained within \cite{BT,HW}.
\end{remark}

\section{`Affine' Setting (3-folds) \`{a} la Keating--Smith and Shimpi}

The setting of this section is $3$-fold flopping contractions $\scrX\to\Spec \scrR$ for which the exceptional curve $\Curve$, given reduced scheme structure, is  isomorphic to $\mathbb{P}^1$. These are known as single-curve flops, and they form a large, in fact uncountable, class of examples, even if we impose the additional assumption that $\scrX$ is smooth.  As in the previous section, consider the category of compactly supported objects
\[
\scrD=\{ a\in\Db(\coh \scrX)\mid \Supp a\subseteq \Curve\}.
\]
The simplest case is that of the Atiyah flop, where $\scrR=\bC\llsq u,v,x,y\rrsq /(uv-xy)$, and $\scrX$ is obtained by blowing up $\Spec\scrR$ at either the ideal $(u,x)$ or $(u,y)$.  

\subsection{Keating--Smith}\label{subsec:KS} The classification of spherical objects in $\scrD$ for the Atiyah flop is due to Keating and Smith \cite{KeatingSmith}. Their paper is remarkable since it was the very first 3-fold `affine' example, which overcame many of the issues with generalising \cite{IU} as explained above, and furthermore their proof technique does not follow the common strategy of \S\ref{sec:CommonStrategy}.

We very briefly overview the strategy in \cite{KeatingSmith} here.  The main idea is to turn everything on its head, and to conclude properties of objects as a consequence of establishing recognition theorems for group elements. The common strategy in \S\ref{sec:CommonStrategy} is, in some sense, the opposite direction to this.

For the Atiyah flop, the category $\scrD$ has some standard spherical objects $\scrS_i=\scrO_{\Curve}(i)$. Now, given an arbitrary spherical object $S\in\scrD$, crucially there is an associated spherical twist $F\coloneq T_S$.  Since this acts trivially on K-theory, the functor $F$ belongs to the subgroup of autoequivalences that we already know \cite[4.11]{KeatingSmith} (namely, spherical twists in the $\scrS_i$, and even shifts) and so the trick is to use this information to deduce that $S$ is a spherical object that we already know.

The proof goes via dynamics.  Via the defining exact triangle of the spherical twist $F=T_S$, it follows that the Floer cohomology $\mathrm{HF}(\scrS_i,F^n\scrS_j; \bZ/2\bZ)$ grows at most linearly with $n$.  Then the main result \cite[6.18]{KeatingSmith}, which uses the Nielsen--Thurston theorem, asserts that up to shift $F=T_S$ is a power of a spherical twist in an element of $\scrS$, which by definition is the set containing all images of all $\scrS_i$ under the action of the group generated by the spherical twists in all the $\scrS_j$. In other words, $T_S$ is, up to shift, a power of a spherical twist that we already know.  

From there, the hard work has been done, and some standard reduction steps (involving e.g.\ tensoring by line bundles) reduces the problem to when $S$ belongs to the null category $\scrC$.  But that situation is known by previous sections, and thus the classification of spherical objects in $\scrD$ follows; the precise statement of the result can be found as the $\ell=1$ case of \ref{thm:Parth} below.

\begin{remark}
From the viewpoint of more general objects, including those considered in \S\ref{subsec:Parth} directly below, the main limitation of the above technique is found in the very first step, where the arbitrary spherical object $S$ gives a group element $T_S$.  Of course, objects do not generate autoequivalences (and thus group elements) in general. Regardless, it is intriguing whether or not there are replacements to the recognition theorem \cite[6.18]{KeatingSmith} in the more general settings below, and this question seems to be intertwined with various open questions on big mapping class groups.
\end{remark}

\subsection{Shimpi}\label{subsec:Parth}
Returning to the general setup of a single-curve flop $\scrX\to\Spec\scrR$, where $\scrX$ may have Gorenstein terminal singularities, recently Shimpi \cite{Parth1, Parth2} has classified spherical (and brick) complexes within the category $\scrD$ of compactly supported objects on $\scrX$.

To set notation and language, to such a flopping contraction is associated Koll\'ar's length invariant $\ell$ \cite{KollarFlops}, which is an integer satisfying $1\leq \ell\leq 6$. Roughly speaking this measures the complexity; the Atiyah flop has $\ell=1$. Given this number, there are natural sheaves 
\[
\scrO_\Curve, \scrO_{2\Curve},\hdots, \scrO_{\ell\Curve},
\]
which encode the successive thickenings of the curve $\Curve\cong\mathbb{P}^1$ within $\scrX$.  When $\ell\geq 5$, by \cite{DW5} there is an additional sheaf, given as the unique non-split extension
\[
0\to\scrO_{3\Curve}\to\scrZ\to\scrO_{2\Curve}\to 0.
\]

\begin{remark}
Here it becomes important to be able to classify more general objects.  Each of the $\scrO_{i\Curve}$ turns out to induce an autoequivalence, after noncommutative deformation \cite{DW5}. But, as explained in \cite[7.5]{DW5}, very few are spherical. Indeed, when e.g.\ $\ell=4$ only $\scrO_{4\Curve}$ is a spherical object.
\end{remark}

With the language set, the following is the main result of \cite{Parth2}, which recovers the classification result (albeit not the results used in the proof) in \S\ref{subsec:KS} as the case $\ell=1$.

\begin{thm}[{\cite{Parth2}}]\label{thm:Parth}
If $x\in \scrD$ is a brick, then up shifts, and iterated applications of mutation functors (which includes line bundle twists) and their inverses, $x$ is either
\begin{enumerate}
\item\label{thm:Parth 1} one of the sheaves $\scrO_\Curve, \scrO_{2\Curve},\hdots, \scrO_{\ell\Curve}$, or possibly $\scrZ$ when $\ell\geq 5$, or
\item\label{thm:Parth 2} the Grothendieck dual of one of the sheaves in \eqref{thm:Parth 1},  or
\item\label{thm:Parth 3} a skyscraper sheaf at a closed point on $\scrX$, or on the flop $\scrX^+$.
\end{enumerate}
\end{thm}
A classification of all objects which satisfy $\Hom_\scrD(x,x[i])=0$ for all $i<0$ also exists \cite[Thm A]{Parth2}, but we refrain from stating that here.

We briefly sketch the strategy behind \ref{thm:Parth}, which is broadly similar to that of the category $\scrC$ outlined in \S\ref{subsec:1.4CLass}, but with significantly greater challenges.  The induction is again on the homological spread with respect to the standard heart induced by the noncommutative resolution.  The key is that now, on this larger category $\scrD$, the induction step becomes much more complicated, in part because the object $x$ must be attacked using both positive atoms `from below' and negative atoms `from above', but also simply because the category $\scrD$ contains more things, so fundamentally there are more options for what can happen.

The very rough overview is that the induction step works since Shimpi understands all intermediate t-structures, given that was the main result of his previous paper \cite{Parth1}.  In fact, the paper \cite{Parth1} works for all $3$-fold flops, and the point is that passing from $\scrC$ to $\scrD$ means replacing the $\tau$-tilting finite and silting-discrete technology of the previous sections (which assert that there are finitely many intermediate t-structures in $\scrC$, all of which are algebraic) by the description of intermediate t-structures in $\scrD$ in \cite{Parth1}. However, this replacement comes at a cost, since $\scrD$ contains many intermediate hearts that are not algebraic, and so care has to be taken if we land in one of these hearts, since in particular it is not so clear how the induction proceeds.  On the positive side, from \cite{Parth1} it turns out that the non-algebraic hearts only come in two basic flavours: straight-up coherent sheaves, and perverse sheaves on partial contractions (known as \emph{mixed hearts}).  In the case of single-curve flops, the mixed hearts do not exist, making the analysis easier, and thus making it possible to establish the induction step. Thus, \ref{thm:Parth} so far works only in the single curve setting.  

\begin{remark}
Given that \cite{Parth1} works in full generality, it seems promising to hope for a more general version of \ref{thm:Parth} which classifies all brick objects for all multi-curve $3$-fold flops, using the same broad strategy.  That said, some of the mixed hearts remain stubbornly problematic, and their existence currently interferes with some of the more subtle support arguments made in \cite{Parth2}.
\end{remark}



\begin{thebibliography}{AMY}
\bibitem[AIR]{AIR}
T.~Adachi, O.~Iyama and I.~Reiten, \emph{$\tau$-tilting theory}, Compos.\ Math.\ \textbf{150} (2014), no.~3, 415--452. 

\bibitem[AMY]{AMY}
T.~Adachi, Y.~Mizuno, and D.~Yang, \emph{Discreteness of silting objects and t-structures in triangulated categories}, Proc.\ Lond.\ Math.\ Soc.\ (3) \textbf{118} (2019), no.~1, 1--42. 

\bibitem[AH]{AH2}
T.~Aihara and T.~Honma, \emph{When is the silting-discreteness inherited?},
Nagoya Math.\ J.\ \textbf{256} (2024), 905--927.

\bibitem[AI]{AI}
T.~Aihara and O.~Iyama, \emph{Silting mutation in triangulated categories}, J.\ Lond.\ Math.\ Soc.\ (2) \textbf{85} (2012), no.~3, 633--668.

\bibitem[AM]{AM}
T.~Aihara and Y.~Mizuno, \emph{Classifying tilting complexes over preprojective algebras of Dynkin type}, Algebra Number Theory \textbf{11} (2017), no.~6, 1287--1315.

\bibitem[AlN]{AlN}
S.~Al-Nofayee, \emph{Simple objects in the heart of a t-structure}, J.\ Pure Appl.\ Algebra \textbf{213} (2009), no.~1, 54--59.

\bibitem[A1]{Asai}
S.~Asai, \emph{Semibricks}, Int.\ Math.\ Res.\ Not.\ IMRN 2020, no.~16, 4993--5054. 

\bibitem[A2]{August1} 
J.~August, \emph{On the finiteness of the derived equivalence classes of some stable endomorphism rings}, Math.\ Z.\ \textbf{296} (2020), no.~3-4, 1157--1183.

\bibitem[A3]{August} 
J.~August, \emph{The Tilting Theory of Contraction Algebras},  Adv.\ Math.\ \textbf{374} (2020), 107372, 56 pp. 

\bibitem[AD]{AD}
J.~August and A.~Dugas, \emph{Silting and tilting for weakly symmetric algebras},  
Algebr.\ Represent.\ Theory \textbf{26} (2023), no.~1, 169--179.

\bibitem[BT]{BT}
C.~Brav and H.~Thomas, \emph{Braid groups and Kleinian singularities}, Math.\ Ann.\ \textbf{351} (2011), no.~4, 1005--1017.

\bibitem[BDL]{BDL}
A.~Bapat, A.~Deopurkar and A.~M.~Licata, \emph{Spherical objects and stability conditions on 2-Calabi--Yau quiver categories}, Math.\ Z.\ \textbf{303} (2023), no.~1, Paper No.~13, 22 pp. 
 
\bibitem[B1]{B02} 
T.~Bridgeland, \emph{Flops and derived categories}, Invent.\ Math.\ \textbf{147} (2002), no.~3, 613--632. 

\bibitem[B2]{Bstab} 
T.~Bridgeland, \emph{Stability conditions on triangulated categories},  Ann.\ of Math.\ (2) \textbf{166} (2007), no.~2, 317--345.

\bibitem[C]{Chen}
J.-C.~Chen, \emph{Flops and equivalences of derived categories for threefolds with only terminal Gorenstein singularities}, J. Differential Geom. \textbf{61} (2002), no. 2, 227--261.

\bibitem[DIJ]{DIJ}
L.~Demonet, O.~Iyama and G.~Jasso, \emph{$\tau$-tilting finite algebras, bricks, and $g$-vectors}, Int.\ Math.\ Res.\ Not.\ IMRN 2019, no.~3, 852--892. 

\bibitem[DWZ]{DWZ}
H.~Derksen, J.~Weyman, and A.~Zelevinsky, \emph{Quivers with potentials and their representations.\ I.\ Mutations.} Selecta Math.\ (N.S.) \textbf{14} (2008), no.~1, 59--119.
 
\bibitem[DW1]{DW1}
W.~Donovan and M.~Wemyss, \emph{Noncommutative deformations and flops}, Duke Math.\ J.\ \textbf{165} (2016), no.~8, 1397--1474. 

\bibitem[DW3]{DW3}
W.~Donovan and M.~Wemyss, \emph{Twists and braids for general $3$-fold flops}, J.\ Eur.\ Math.\ Soc.\ (JEMS), \textbf{21} (2019), no.~6, 1641--1701.

\bibitem[DW5]{DW5}
W.~Donovan and M.~Wemyss, \emph{Stringy K\"ahler moduli, mutation and monodromy}, J.\ Differential Geom.\ \textbf{129} (2025), no.~1, 115--164.

\bibitem[EJR]{EJR}
F.~Eisele, G.~Janssens and T.~Raedschelders, \emph{A reduction theorem for $\tau$-rigid modules}, Math.\ Z.\ \textbf{290} (2018), no.~3-4, 1377--1413.

\bibitem[HW1]{HaraWemyss}
W.~Hara and M.~Wemyss, \emph{Spherical objects in dimensions two and three}, to appear J.\ Eur.\ Math.\ Soc.\ (JEMS).

 \bibitem[HW2]{HW} 
 Y.~Hirano and M.~Wemyss, \emph{Faithful actions from hyperplane arrangements}, Geom.\  Topol.\ \textbf{22} (2018), no.~6, 3395--3433. 

\bibitem[H]{Huy} D.~Huybrechts, \emph{Stability conditions via spherical objects}, Math.\ Z.\ \textbf{271} (2012), no. 3--4, 1253--1270.

\bibitem[IU]{IU}
A.~Ishii, and H.~Uehara, \emph{Autoequivalences of derived categories on the minimal resolutions of $A_n$-singularities on surfaces}, J.\ Differential Geom.\ \textbf{71} (2005), no.~3, 385--435.

\bibitem[IUU]{IUU}
A.~Ishii, K.~Ueda and H.~Uehara, \emph{Stability conditions on $A_n$-singularities}, J.\ Differential Geom.\ \textbf{84} (2010), no.~1, 87--126. 

\bibitem[IW1]{IW1}  
 O.~Iyama and M.~Wemyss, \emph{Maximal modifications and Auslander-Reiten duality for non-isolated singularities}. Invent.\ Math.\ \textbf{197} (2014), no.~3, 521--586.
 
\bibitem[IW2]{IW9}
O.~Iyama and M.~Wemyss, \emph{Tits cones intersections and applications}, \href{https://www.maths.gla.ac.uk/~mwemyss/MainFile_for_web.pdf}{\texttt{preprint}}. 

\bibitem[K1]{Kawatani}
K.~Kawatani,
\textit{Pure sheaves and Kleinian singularities},
Manuscripta Math. \textbf{160} (2019), no. 1--2, 65--78.

\bibitem[KS]{KeatingSmith}
A.~Keating and I.~Smith, \emph{Symplectomorphisms and spherical objects in the conifold smoothing}, Compos.\ Math.\ \textbf{160} (2024), no.~11, 2738--2773.

\bibitem[K2]{Kimura}
Y~Kimura, \emph{Tilting and silting theory of Noetherian algebras},
Int.\ Math.\ Res.\ Not.\ IMRN 2024, no.~2, 1685--1732.

\bibitem[KY]{KY}
S.~Koenig and D.~Yang, \emph{Silting objects, simple-minded collections, t-structures and co-t-structures for finite-dimensional algebras}, Doc.\ Math.\ \textbf{19} (2014), 403--438. 

\bibitem[K3]{KollarFlops}
J.~Koll\'ar, \emph{Flops}. Nagoya Math.\ J.\ \textbf{113} (1989), 15--36.


\bibitem[KM]{KM}
J.~Koll\'ar and S.~Mori,  \emph{Birational geometry of algebraic varieties}, Cambridge Tracts in Math., \textbf{134}, Cambridge Univ. Press, Cambridge, 1998. 


\bibitem[M]{Mizuno}
Y.~Mizuno,
\emph{Classifying $\tau$-tilting modules over preprojective algebras of Dynkin type},
Math.\ Z.\ \textbf{277} (2014), no. 3--4, 665--690.

\bibitem[PSZ]{David}
D.~Pauksztello, M.~Saor\'{i}n, A.~Zvonareva, \emph{Contractibility of the stability manifold for silting-discrete algebras}, Forum Math.\ \textbf{30} (2018), no.~5, 1255--1263. 

\bibitem[P]{Pinkham}
H.~Pinkham, \emph{Factorization of birational maps in dimension 3}, Singularities (P. Orlik, ed.), Proc. Symp. Pure Math., vol. 40, Part 2, 343--371, American Mathematical Society, Providence, 1983.

\bibitem[R]{Rickard}
J.~Rickard, \emph{Equivalences of derived categories for symmetric algebras}, J.\ Algebra \textbf{257} (2002), no.~2, 460--481.

\bibitem[ST]{ST}
P.~Seidel and R.~Thomas, \emph{Braid group actions on derived categories of coherent sheaves}, Duke Math.\ J.\ \textbf{108} (2001), no.~1, 37--108.

\bibitem[S1]{Parth1}
P.~Shimpi, \emph{Torsion pairs and 3-fold flops}, \href{https://arxiv.org/abs/2502.05146}{\texttt{arXiv:2502.05146}}.

\bibitem[S2]{Parth2}
P.~Shimpi, \emph{Simple complexes on a flopped curve}, \href{https://arxiv.org/abs/2508.06437}{\texttt{arXiv:2508.06437}}.

\bibitem[SW]{SW}
I.~Smith and M.~Wemyss, \emph{Double bubble plumbings and two-curve flops}, Selecta Math.\ (N.S.) \textbf{29} (2023), no.~2, Paper No.~29, 62 pp. 

\bibitem[VdB]{VdB}
M.~Van den Bergh, \emph{Three-dimensional flops and noncommutative rings}, 
Duke Math.\ J.\ \textbf{122}  (2004), no.~3, 423--455.

\bibitem[W]{HomMMP}
M.~Wemyss, \emph{Flops and Clusters in the Homological Minimal Model Program}, Invent.\ Math.\ \textbf{211} (2018), no.~2, 435--521.

\bibitem[Z]{Hao}
H.~Zhang, \emph{Local forms for the double $A_n$ quiver}, \href{https://arxiv.org/abs/2412.10042}{\texttt{arXiv:2412.10042}}.







\end{thebibliography}
\end{document}